\documentclass[nonblindrev]{informs3}

\DoubleSpacedXI 

\usepackage{xr}
\usepackage{algorithm}
\usepackage{algpseudocode}
\usepackage{endnotes}
\usepackage{caption}
\usepackage{subcaption}
\let\footnote=\endnote
\let\enotesize=\normalsize
\def\notesname{Endnotes}%
\def\makeenmark{\hbox to1.275em{\theenmark.\enskip\hss}}
\def\enoteformat{\rightskip0pt\leftskip0pt\parindent=1.275em
	\leavevmode\llap{\makeenmark}}

\usepackage{natbib}
\bibpunct[, ]{(}{)}{, }{a}{}{, }%
\def\bibfont{\small}%
\newenvironment{assumptionbis}[1]
{\renewcommand{\theassumption}{\ref{#1}$'$}%
	\addtocounter{assumption}{-1}%
	\begin{assumption}}
	{\end{assumption}}

\TheoremsNumberedThrough     
\ECRepeatTheorems

\EquationsNumberedThrough    

\begin{document}
	
	
	\RUNAUTHOR{Jin et al.}
	
	\RUNTITLE{CONE for prescriptive analysis}
	
	\TITLE{Contextual Optimizer through Neighborhood Estimation for prescriptive analysis}
	
	\ARTICLEAUTHORS{%
		\AUTHOR{Xiao Jin}
		\AFF{SIA-NUS Digital Aviation Corporate Laboratory, National University of Singapore, Singapore, 117602, SINGAPORE, \\ \EMAIL{isejinx@nus.edu.sg}} 
	\AUTHOR{Yichi Shen}
	\AFF{C3 AI, Singapore, 068877 \EMAIL{isesy@nus.edu.sg}}
	
	\AUTHOR{Loo Hay Lee}
	\AFF{Industrial Systems Engineering and Management, 
		National University of Singapore, 
		Singapore, 117576, SINGAPORE, \EMAIL{iseleelh@nus.edu.sg}}
	\AUTHOR{Christine A. Shoemaker}
	\AFF{Industrial Systems Engineering and Management, 
		National University of Singapore, 
		Singapore, 117576, SINGAPORE, \EMAIL{shoemaker@nus.edu.sg}}
} 

\ABSTRACT{%
	We address the challenges posed by heteroscedastic noise in contextual decision-making. We propose a consistent Shrinking Neighborhood Estimation (SNE) technique that successfully estimates contextual performance under unpredictable variances. Furthermore, we propose a Rate-Efficient Sampling rule designed to enhance the performance of the SNE. The effectiveness of the combined solution ``Contextual Optimizer through Neighborhood Estimation"(CONE) is validated through theorems and numerical benchmarking. The methodologies have been further deployed to address a staffing challenge in a hospital call center, exemplifying their substantial impact and practical utility in real-world scenarios.
}%


 \HISTORY{This paper was
	first submitted on Aug 3, 2023.}

\maketitle

%


\section{Introduction}
\label{section:intro}

The inherent contextuality of operations problems establishes that an optimal decision $x$, while ideal in one context $y$, may not universally apply. This paradigm of ``context-dependency" demands an optimal rule mapping each context to an associated optimal decision (e.g., $x^*(y)$) based on the contextual performance $f_x(y)$. Prescriptive analysis  \citep{bertsimas2020predictive} can be leveraged when ample data is available. The concept involves employing statistical learning to train a proxy $\hat{x}(y)$ from labeled performance observations $\tilde{G}_x(y)$, estimating the conditional optimal decision $x^*(y)$ for each $y$ context. However, developing a robust proxy necessitates substantial performance observations spanning a diverse range of decision and context combinations. This challenge intensifies considering the inherent heteroscedastic noise accompanying observations $\tilde{G}_x(y)$, irrespective of their derivation from historical data, field experiments, or computer simulations. In this research, we concentrate on  \textit{the generation of data from a heteroscedastic simulation oracle to expedite the proxy training process}, thereby advancing the capabilities of prescriptive analysis in operations research.

The employment of simulation to enhance prescriptive analysis, as suggested by \cite{hong2019offline}, seeks to optimize context-dependent performance through strategic simulation experiments. However, this approach presents substantial challenges. Firstly, achieving a consistent performance estimation across the entire spectrum of decision and context combinations is challenging. While the Sample Average Approximation (SAA) on replicated $(x, y)$ data can provide reliable estimation in limited combinations, it struggles in cases where these combinations are numerous or infinite, resulting in a limited coverage of the whole decision-context space. Secondly, the development of a budget allocation strategy, even if rooted in a well-designed performance estimator, is difficult due to the high-dimensionality of the decision-context joint space and the presence of heteroscedastic noise.

In this study, we address the intricate nature of contextual problems characterized by discrete, categorical decisions with either discrete or continuous contexts. Such problem settings, common in areas like inventory management and staff scheduling (contingent on demand prediction), often cloud the connections between decisions and conditional outcomes due to their inherent complexity, infinite context possibilities, and notable heteroscedasticity. Given these challenges, our research focuses on: 1. formulating an estimator proficient in addressing heteroscedasticity and utilizing single observations for each decision-context pair, and 2. crafting a budget allocation strategy based on this performance estimator to minimize overall decision.
\subsection{Challenges from infinite context and heteroscedastic noise}
\begin{figure}[!h]
	\centering
	\includegraphics[width=0.7\textwidth]{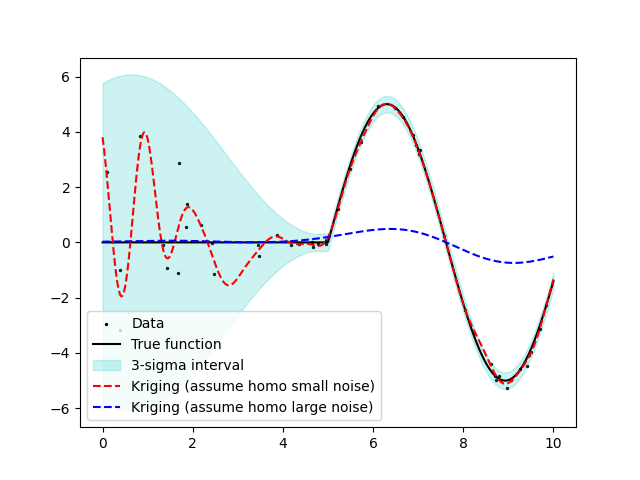}
	\caption{Illustration of Kriging model limitations in single observation heteroscedastic data interpolation. 
		Data points are generated by introducing heteroscedastic noise to the true function (black line). Interpolations (red and blue dashed lines) are carried out using the Kriging model under the assumption of different homoscedastic noise. Due to significant noise, the surrogate by assuming small noise on the left part of the space inaccurately captures the true function, leading to skewed estimations. Conversely, the right part of the space features minimal noise, and assuming high levels of homoscedastic noise leads to distortion in estimations.}\label{Figure: challenge}
\end{figure}

To estimate the conditional performance on infinite context space, increasing the data coverage on decision-context space is crucial, favoring single simulation assigned for each decision-context combination. However, due to the existence of heteroscedastic noise, it is also important to do simulation run at or near observed decision-context to mitigate the estimation noise. These leads us to three common strategies in solving this trade-off problem. Firstly, one can simply ignore the heteroscedastic noise and predict performance by assuming homoscedastic noise. However, this strategy has significant flaws.  As shown in Figure \ref{Figure: challenge}, using homoscedastic assumptions on heteroscedastic data, such as with Kriging, will invariably result in imprecise estimations on part of the domain. Secondly, we could also estimate heteroscedastic noise by taking average among replications at the same decision-context point. However, multiple observations per decision-context point are compulsory for this strategy, limiting exploration on a vast high-dimensional context space.

This leads us to the third strategy: neighborhood estimation. It uses observations from close enough contexts to predict the performance of the unobserved one. 
On one hand, neighborhood estimation has an intrinsic variance estimator provided by the sample variance of the neighborhood samples, enabling allocation rule that utilize information of the local noise level. On the other hand, it, at the same time, allows single observation which provides maximum flexibility to maximize data coverage on decision-context space.

Designing a good neighborhood is also a hurdle. Since the observations in the neighborhood are used, on top of the heteroscedastic noise, it introduces bias from the neighborhood. The well-known k-NN regression\citep{mack1981local} uses simple/weighted average among k nearest observations to estimate the current point. By increasing the observations in the neighborhood, the size of the neighborhood region will decrease, mitigating the bias. However, k-NN regression cannot provide a consistent estimation in noisy cases since the size of neighborhood samples is always fixed to be $k$ and the noise error cannot be averaged out. For neighborhood estimation that can achieve consistency, shrinking ball (\cite{linz2017computational},\cite{kiatsupaibul2018single},\cite{kiatsupaibul2020single}) controls both the heteroscedastic noise and the bias by adjusting the size of the neighborhood so that the number of observations within goes to infinity while the size of the neighborhood goes to zero.

In this study, we construct an innovative neighborhood rule founded on three fundamental principles: 1) computational simplicity of estimators; 2) consistent estimation for both expected performance and variance of heteroscedastic noise; 3) minimal data prerequisites, such as multiple observations per decision-context combination. Considering the categorical nature of decisions and the continuous character of contexts within our problem, we designate a neighborhood for each decision-context pair as the set sharing an identical decision and a similar context.
\subsection{Budget allocation for prescriptive analysis}
The goal of budget allocation for prescriptive analysis is to optimize the simulation data generation so that the decision loss caused by implementing the proxy trained by the data is minimized. There is few literature that address the exact general problem we raised. Still, we review two streams of literature associated with conventional budget allocation problem in simulation-based optimization to help locate our position in the literature. The first stream is the Multi-Arm Bandit (MAB) and the second is Ranking and Selection (R\&S). 

The goal of MAB is to maximize the accumulated reward in a series of repeated decisions by choosing one decision (arm) each round \citep{bouneffouf2019survey}. State-dependent/contextual MAB problem is also studied  \citep{lu2010contextual}, where at each round, the state/side-information is revealed to benefit the decision-making. From the MAB perspective, our study focuses on the pure exploration of continuous state-dependent MAB. This means the objective is to find the optimal response $x^*(y)$ as exactly as possible regardless of the reward accumulated on the way of learning. It is similar to the pure exploration MAB  \citep{bubeck2011pure} or best arm identification problem \citep{audibert2010best}  where the goal is to find the best decision as well. However, little literature has discussed the state-dependent best arm identification MAB with heteroscedastic noise.

The goal of R\&S is to find the optimal decision among finite options\citep{hong2021review}. There are fundamentally two types of methods: fix precision and fix budget. Fix precision methods such as the Indifference Zone (IZ) \citep{kim2001fully} intend to guarantee a frequentist probability of selecting the real best decision. Fix budget methods such as Expected Value of Information(EVI) \citep{chick2010sequential} and Knowledge Gradient (KG) \citep{frazier2008knowledge} intend to allocate each additional observation to benefit the posterior probability of selecting the optimal myopically. Another popular fix-budget method is the Optimal Computing Budget Allocation (OCBA) \citep{chen2011stochastic} which uses the fixed allocation ratio that asymptotically optimizes the Large Deviation (LD) \citep{glynn2004large} rate as a guide to allocate budget. Based on this categorization, the budget allocation derived in this paper is a sequential fixed budget algorithm.

Although contextual Ranking and Selection (R\&S) has garnered attention recently, few research addressed problems in continuous context space. In  \cite{pearce2018continuous}, a Bayesian framework was proposed and Kriging is adopted to describe the posterior distribution of performance on unobserved $(x,y)$ pairs. In their work, they consider the case where decisions are continuous while scenarios can be discrete or continuous. The proposed myopic budget allocation maximizes overall Expected Improvement (EI). Later in \cite{li2020context}, a similar Bayesian framework is used and the states and decisions are both finite. A dynamic programming perspective is introduced to define the optimal budget allocation rule. Their goal is to maximize the minimum Probability of Correct Selection among the states. The most recent study (\cite{cakmak2022contextual}, \cite{ding2022knowledge}) both use Gaussian process regression to anticipate the performance. \cite{cakmak2022contextual} allocates budget according to optimize the decay rate of the PFS and the consistency and the decay rate of the algorithm are proved. \cite{ding2022knowledge} uses Knowledge Gradient methods to maximize the benefits of each budget allocated.  All these Gaussian/Kriging regression methods face the same challenge illustrated in Figure \ref{Figure: challenge} when dealing with heteroscedastic noise. All these methods either assume that the noise is homogeneous or heteroscedastic with known variance. If the variance of the heteroscedastic noise is unknown, replications at the same point are still inevitable. \cite{gao2019selecting} and \cite{jin2019optimal} intend to maximize the asymptotic decay rate of the overall false identification of the conditional optimal decisions. These methods can guarantee the optimal decay rate of selecting false decisions averaged on states, however, are only applied where state space is finite and small. In addition, multiple replications at the same $(x, y)$ are still inevitable.

For research addressed continuous context space, \cite{shen2021ranking} used linear regression to model context-performance relationships, offering an algorithm that can ensure overall decision loss, termed Probability of False Selection (PFS).  Despite considering a continuous context space, this method needs replications on a fixed set of states to facilitate linear regression against heteroscedastic noise. This puts the validity of the method heavily on the degree of linearity of the context-performance relation. Furthermore, the approach may allocate the budget conservatively, lacking adaptiveness if large noise levels are detected in any decision-context pair.

\subsection{Gaps and Contributions}
The existing literature largely treats simulation noise as either homoscedastic or heteroscedastic with known variance. A limited methods does consider the unknown variance of heteroscedastic noise, but they are typically constrained by small, finite context spaces or rely heavily on specific structural assumptions, such as linearity. The condition of unknown heteroscedastic noise variance almost always leads to repeated observations at identical decision-context pairs. Currently, there is a marked absence of methods capable of managing high heteroscedasticity with a single observation per decision-context pair. 

Building upon the identified research gaps, this paper brings forth three significant contributions to the field:
\begin{enumerate}
	\item We develop a novel neighborhood estimation method that is adaptive to heteroscedastic noise, liberating the need for replications at the same $(x, y)$ combination. We demonstrate the consistent nature of this estimation, both in terms of expected performance and variance, ultimately facilitating the accurate recovery of real conditional optimal decisions for any given state $y$.
	\item We introduce an efficient sequential rule for budget allocation, underpinned by proven optimal convergence rate of the decision loss. This rule is adept at directing the budget towards those $(x, y)$ combinations that are most pivotal to overall decision loss, optimizing resource distribution. Comparative evaluations against established algorithms further underscore the superior performance of our proposed method.
	\item The utility and applicability of our method are demonstrated through practical implementation in a hospital call center staffing problem. The objective is to prescribe the best staffing allocation based on caller arrival and patience profile, revealing the method's potential in delivering high-quality prescriptive analysis even in high-dimensional, heteroscedastic scenarios.
\end{enumerate}

The structure of this paper is as follows: In Section \ref{section:problem}, we present a comprehensive problem definition, highlighting two core challenges: optimality inference and budget allocation. We also delineate the technical assumptions underlying our analysis. In Sections \ref{section:methodology_SNE} and \ref{section:methodology_Allocation}, we address these challenges in detail. Section \ref{section:methodology_SNE} outlines the specifics of the Shrinking Neighborhood Estimation (SNE) methodology, including a rigorous convergence analysis. The subsequent Section \ref{section:methodology_Allocation} introduces CONE, our novel budget allocation strategy, derived from an optimization on the convergence rate. Section \ref{section: numericals} positions our proposed budget allocation strategy against benchmark rules through numerical comparisons, substantiated by a theoretical case and a practical application in a hospital call center. These case studies underscore the substantial benefits of integrating our new estimation method with the budget allocation rule. The paper concludes in Section \ref{section: conclusion}, where we revisit the main contributions of our research, reinforcing its impact and potential implications for the field.

\section{Problem Definition}
\label{section:problem}
Suppose a decision maker can choose decision $x$ among $\chi$ where $|\chi|$ is finite, in response to context $y\in Y$ where $Y\subset\mathbb{R}^{d}$ is bounded. For each decision $x$, the expected performance on context $y$ is a real-valued function $f_x: Y\mapsto \mathbb{R}$. The decision maker aims at achieving a perfect proxy, which assigns each context with the optimal decision:
\begin{equation}
	\label{equation:simopt}
	x^*(y)=\arg\min_{x\in \chi}f_{x}(y)\quad\forall y\in Y,
\end{equation}  
However, $f_x(y)$ is unknown, and can only be estimated via stochastic simulation observations. Each simulation observation $\gamma:=((x, y), \tilde{G}_x(y))$ consists of a the specified $(x, y)$ combination and one realization $\tilde{G}_x(y):=f_x(y)+\epsilon_x(y)$, where $\epsilon_x(y)$ is the random error.  We use $``\cdot_{[i]}"$ to label each run in sequence. Therefore, for the $i^{th}$ run $\gamma_{[i]}$, the simulation run is $\gamma_{[i]}:=((x_{[i]}, y_{[i]}), \tilde{G}_{[i]})$, where $ \tilde{G}_{[i]}={f}_{[i]}+{\epsilon}_{[i]}$ and ${f}_{[i]}:=f_{x_{[i]}}(y_{[i]})$. The dataset $\mathcal{D}$ is a sequence of labeled experiments result (i.e. $\mathcal{D}:=\big(\gamma_{[1]}, \gamma_{[2]}, \dots, \gamma_{[N]}\big)$). Given $\mathcal{D}$, if $\hat{f}_x(y)$ is specified, the inferred conditional optimal decision is:
\begin{equation}
	\hat{x}(y):=\arg\min_{x\in \chi} \hat{f}_x(y)\quad \forall y\in Y.\label{subopt}
\end{equation}We force uniqueness of $\hat{x}(y)$ by randomly select one $x\in \hat{x}(y)$ if two or more $x$ share the same optimal $\hat{f}_x(y)$ on one $y$. 

To quantify the performance of $\hat{x}(y)$, a total decision loss is defined. Let $f^*(y):=f_{x^*(y)}(y)$ and $f^\dagger(y):=f_{\hat{x}(y)}(y)$. Conditioned on $y$, a decision loss incurs if $f^\dagger(y)\neq {f^*(y)}$. For a general metric function $d(\cdot,\cdot)$, the loss at $y$ is defined as $d(f^\dagger(y),{f^*(y)})$. The expected decision loss of using $\hat{x}(y)$ conditioned on $y$ is
\begin{equation}
	\upsilon(y):=E_\mathcal{D} d(f^\dagger(y),{f^*(y)}),\forall y\in Y.
\end{equation}
Our definition of decision loss is based on a general metric function. Different specifications of the metric functions are found in existing literature. If $d(f^\dagger(y),{f^*(y)})=1$ when $f^\dagger(y)\neq{f^*(y)}$ and $d(f^\dagger(y),{f^*(y)})=0$ when $f^\dagger(y)={f^*(y)}$, $\upsilon$ is called to ``Probability of False Selection" (PFS) \citep{wu2018analyzing}. Alternatively, if $d(f^\dagger(y),{f^*(y)})=|f^\dagger(y)-{f^*(y)}|$, $\upsilon$ is equivalent to the Expected Opportunity Cost (EOC) \citep{branke2007selecting}. The total decision loss of implementing $\hat{x}(y)$ is defined as the integral of $\upsilon(y)$ on context space. To introduce weight on different contexts, one can equipped $Y$ with a probability distribution $P_Y$. $P_Y$ represents the practitioner's belief in the occurrence likelihood of events. Therefore, the total decision loss caused by implementing $\hat{x}(y)$ is:
\begin{equation}
	\varpi=\int_{y\in Y} \upsilon(y) dP_Y.
\end{equation}\label{definition: varpi}

\subsection{Performance estimation and budget allocation}
There are two identified problems in achieving a perfect proxy. First is how to estimate the ${f}_x(y)$ based on a given $\mathcal{D}$. Second is how to generate $\mathcal{D}$ knowing that $\mathcal{D}$ will be applied to estimate ${f}_x(y)$.

For the performance estimation problem, let $\hat{f}_x(y)$ denote an estimator of  $f_x(y)$ at $(x,y)$. Notably, given the benefits of neighborhood estimation, $\hat{f}_x(y)$ in this study is determined by averaging neighborhood samples. Constructing a neighborhood structure can be approached in multiple ways.  On one end, maintaining a constant neighborhood region size means that as more observations amass on $\chi\times Y$, the growing neighborhood samples reduce estimation noise, yet the bias remains. Conversely, keeping neighborhood sample size static reduces bias but retains noise. Our tailored neighborhood structure, balancing bias and noise, is elaborated in section \ref{section:methodology_SNE}.

For budget allocation problem, in the case of a simulation study, we are in full control of how data $\mathcal{D}$ is generated. Therefore, the goal of budget allocation is to designed a data generation rule to minimize $\varpi$. In this study, we aim to devise an efficient sequential budget allocation rule to minimize $\varpi$. Under this framework, $\mathcal{D}$ is sequentially updated with each addition of $\gamma_{[i]}$ until a budget $T$ is exhausted. The final $\mathcal{D}$ is then used to produce $\hat{x}(y)$. We employ the superscript $\cdot^{(t)}$ to denote each stage, such that at stage $t$, the budget allocation determines $(x_{[t]}, y_{[t]})$ based on the preceding $\mathcal{D}^{(t-1)}$. Once $(x_{[t]}, y_{[t]})$ is specified, the simulator produces the associated $\tilde{G}_{[t]}$. Subsequently, $\mathcal{D}^{(t)}$ is updated with $\gamma_{[t]}$. The resultant $\hat{x}^{(T)}(y)$ is then implemented in operations. Our goal is to find a budget allocation rule that can minimize $\varpi^{(T)}$.

To facilitate the analysis in the following sessions, we impose the following three assumptions:
\begin{assumption}\label{assumption: lipschitz}
	$\forall x\in\chi$, $f_{x}(y)$ is a uniformly bounded, Lipschitz continuous with constant $L_1>0$ on $Y$.
\end{assumption}
\begin{assumption}
	\label{assumption: variance}
	$\epsilon_x(y)$ follows normal distribution:
	\begin{equation}
		\epsilon_x(y)\sim\mathcal{N}(0, {\sigma_{x}(y)}), 
	\end{equation}
	where $\forall x, y$, ${\sigma_{x}(y)}\in(0, \infty)$. In addition, $\forall x\in\chi$, $\sigma_{x}:Y\mapsto \mathbb{R}$ is Lipschitz continuous with constant $L_2>0$ on $Y$.
\end{assumption}	
\begin{assumption}
	${\epsilon}_{[i]}$ are mutually independent.\label{assumption: independent}
\end{assumption} 
\paragraph{Remark:} For Assumption \ref{assumption: lipschitz}, the uniform boundedness and the Lipschitz continuity enable a consistent neighborhood estimator. For Assumption \ref{assumption: variance}, we assume that the simulation oracle is unbiased. Also, the normality assumption is common in simulation-based optimization research. For non-normal noise, simulations can be batched, with each batch's sample mean representing an observation. Leveraging the Central Limit Theorem, the batch mean will asymptotically approximate a normal distribution as batch size increases. Additionally, the presumption of Lipschitz continuity in the variance function ensures the feasibility of a consistent neighborhood estimation of $\sigma_x(y)$. For Assumption \ref{assumption: independent}, it holds by running each simulation independently.

\section{Shrinking Neighborhood Estimation}
\label{section:methodology_SNE}
In this session, we present a method of establishing $\hat{f}_x(y)$ from a given dataset $\mathcal{D}$, namely Shrinking Neighborhood Estimation (SNE). The basic idea is to estimate $f_x(y)$ at any $(x,y)$ using sample average of $\tilde{G}_{[i]}$ from its ``neighborhood". This neighborhood is defined as observations that has the same $x$ (e.g. $x_{[i]}=x$) and nearby $y$ (e.g. $\|y_{[i]}-y\|<r$ for some radius $r>0$). Due to taking average among neighborhood samples, it introduces two estimation errors. On one hand, $f_{[i]}$ from the neighborhood differs from $f_x(y)$, introducing estimation bias to the sample average. On the other hand, $\tilde{G}_{[i]}$ is noisy itself, thus the sample average is noisy. The challenge of balancing between estimation bias and noise motivates our adaptive design to specify the size of the neighborhood. One additional note is that we will be focusing on analyzing $y$ that belongs to the interior of $Y$ (denoted as $Y^\circ$) because a neighborhood is not always well-defined at the boundary.

On stage $t$, for each $x$, the local sample set of context $y\in Y^\circ$ within radius of $r$ is denoted as $b^{(t)}(x, y;r)$:
\begin{equation}
	b^{(t)}(x, y;r ):=\{\gamma_{[i]}\in\mathcal{D}^{(t)}: x_{[i]}=x,\|y_{[i]}-y\|<r,y_{[i]}\in Y^\circ\}.
\end{equation}
The shrinking process is established based on a sequence of neighborhood, and is specified by a monotone decreasing series $r_m$, $m=1,2,\dots$ defined as:
\begin{equation}
	r_m= c m^{  -\xi  /d},\label{definition: decreasing sequence}
\end{equation}
for some constant $\xi>0$ and $c>0$, of which the technical settings are discussed in the supplement. By definition of $r_m$,  
\begin{equation}
	b_m^{(t)}(x, y):=b^{(t)}(x, y;r_m)
\end{equation}\label{definition: neighborhood samples}specifies a series of nested $m$-labeled sample set (e.g. for $m>m'$, $b_m^{(t)}(x,y)\subset b_{m'}^{(t)}(x,y)$), any of which can be a neighborhood sample set. We call this $b_m^{(t)}(x,y)$ the $m$-set for $x,y$ on stage $t$. For a fixed $t$ and $(x,y)$, the larger the $m$, the smaller the radius $r_m$ and the number of samples within $m$-set are. Intuitively, to guarantee consistent decreasing noise and bias as more observations are accumulated, $m$-set with a larger $m$ is preferred to the extend that the number of samples inside is sufficiently large. To achieve this, we define that $m$-set at $(x,y)$ on stage $t$ is ``activated" if $|b_m^{(t)}(x, y )|\geq m$ so that any $m$-set with at least $m$ observations within are qualified to be a candidate neighborhood set. Notice that $m$-set at $(x,y)$ will be activated at most once throughout the accumulation of observations. Once it is activated on stage $t$, it will always been activated afterwards. The problem left is to determine which $m$-set should be chosen as the neighborhood set for each $(x,y)$ on stage $t$. 
\begin{figure}[!h]
	\centering
	\includegraphics[width=0.8\textwidth]{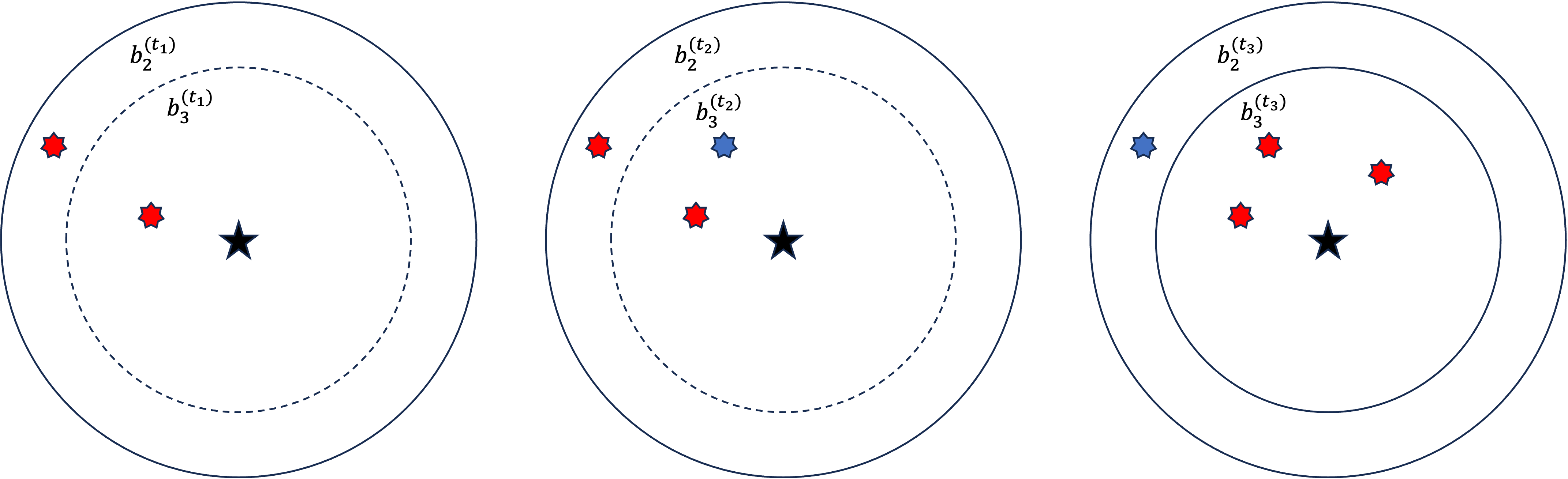}
	\caption{Images showing the neighborhood samples for a given $x$ and $y$ where $Y\subset\mathbb{R}^2$. The black star is $y$. Red stars are neighborhood samples of $(x,y)$ on stage $t$. Blue stars are samples near $y$ that share the same $x$ but not belongs to the neighborhood samples by the SNE rule. From left to right, three stages follow $t_1<t_2<t_3$. 1. The left image shows the case on stage $t_1$ where exactly two observations fall into $b_2$, activating $b_2$, making the neighborhood set $b_{x,y}^{(t_1)}=b_2^{(t_1)}$. 2. The middle image illustrates an intermediary state $t_2$, where the sample set remains unchanged as $b_3$ has yet to include $3$ observations. This leads to $b_{x,y}^{(t_2)}=b_2^{(t_1)}$. 3. The right image shows the final stage $t_3$ when $b_3$ is activated, at which point observations outside $b_3$ are discarded, and $b_{x,y}^{(t_3)}=b_3^{(t_3)}$.}\label{Fig: SNE progression}
\end{figure}

To specify an $m$-set as the neighborhood sample set, we implement a ``passive greedy" rule. On one hand, the rule is greedy in the sense that, for each $(x,y)$, we choose the activated $m$-set with the largest $m$ as the neighborhood samples. On the other hand, the rule is passive in the sense that once a new $m$-set is activated at some $(x,y)$ on stage $t$, the neighborhood samples will remain the same until an $m$-set with a larger $m$ is activated. For each $m$, on stage $t$, denote the first time $t'$ such that $b_m^{(t')}(x,y)$ is activated as $\tau_m^{(t)}(x,y)$:
\begin{equation}
	\tau_m^{(t)}(x,y)=\begin{cases}
		\min\{t'<t:|b_m^{(t')}(x,y)|=m\} & |b_m^{(t)}(x,y)|\geq m\\
		\infty& otherwise
	\end{cases}
\end{equation} At time $t$, denote the last stage with a new activation of an $m$-set as:
\begin{equation}
	\tau^{(t)}_{x,y}:=\max_{m}\{	\tau_m^{(t)}(x,y)<\infty\}.
\end{equation} 
The corresponding $m$ is denoted as
\begin{equation}
	\Psi^{(t)}_{x,y}:=\arg\max_{m}\{	\tau_m^{(t)}(x,y)<\infty\}.
\end{equation} Therefore, given $\mathcal{D}^{(t)}$ at $t$, the neighborhood samples at $(x,y)$ is defined as
\begin{equation}
	b_{x,y}^{(t)}:=b_{\Psi^{(t)}_{x,y}}^{(\tau^{(t)}_{x,y})}(x,y). \label{definition: neighborhood sample}
\end{equation} 
The reason for adopting this rule is twofold. On one hand, the greedy rule will ensure the fastest decrease of the bias by selecting the neighborhood with the smallest radius. On the other hand, in conjunction with the condition that only one observation is acquired on each stage, the passive rule guarantees that the size of the neighborhood samples always equals to the $m$ of the chosen $m$-set, simplifying the convergence analysis. For a comprehensive illustration of this neighborhood rule, please refer to Figure \ref{Fig: SNE progression}.

Based on the definition of $b_{x,y}^{(t)}$, the estimation of the conditional performance, namely ``Shrinking Neighborhood Estimation" (SNE), is the sample average of the neighborhood samples:
\begin{equation}
	\hat{f}_{x}^{(t)}(y):=\frac{1}{\Psi^{(t)}_{x,y}}\sum_{x_{[i]}=x, y_{[i]}\in b_{x, y}^{(t)}} \tilde{G}_{[i]}.
\end{equation}
We also use the sample variance in the neighborhood to estimate the variance:
\begin{equation}
	\hat{\sigma}_{x}^{(t)}(y):=\sqrt{\frac{1}{\Psi^{(t)}_{x,y}-1}\sum_{ x_{[i]}=x,y_{[i]}\in b_{x, y}^{(t)}} [\tilde{G}_{[i]}-\hat{f}_x^{(t)}(y)]^2}.
\end{equation}
The conditional optimal $\hat{x}^{(t)}(y)$ is then inferred from $\hat{f}_{x}^{(t)}(y)$ as suggested in formulation (\ref{subopt}).

There are three advantages of SNE. The first one is computational. Given a new $y$, calculating $\hat{x}^{(t)}(y)$ takes two steps. In step one, the neighborhood sample set $b_{x,y}^{(t)}$ is determined by screening through each observation $\gamma_{[i]}\in \mathcal{D}^{(t)}$ and calculating its distance to the given $y$. Calculation of distance between two points is proportional to the dimension $d$. The time complexity of calculating distance for all samples is $O(td)$. $\hat{f}_x(y)$ is derived by taking average among $\tilde{G}_{[i]}$ within $b_{x,y}^{(t)}$ which is $O(t)$ in the worst case. In step two, $\hat{f}_x(y)$ is sorted among $x$ to indicate $\hat{x}(y)$, which typically costs $O(|\chi|\ln(|\chi|))$. Compared to other inference methods such as Kriging, in step one, they typically requires the inversion of a covariance matrix with the size proportional to $t$ to derive $\hat{f}_x(y)$. The conventional method of matrix inverse typically costs $O(t^3)$. In this regard, our SNE limits the computational effort to compute $\hat{x}(y)$ even when the number of accumulated observations is large. The second advantage is the capability to estimate both the expected performance and the variance of the noise. The estimation of variance can highlight decision-context combinations of which performance requires more observations to estimate. The rule also adaptively specifies the neighborhood for each $x,y$ to guarantee decreasing of the neighborhood size and increasing of the sample size so that both bias and noise of the estimation error keep vanishing. The last advantage is the minimum requirements for data. SNE only requires that for each decision $x$, at least two $y$ are observed to enable local variance estimation. Replications at the same $(x,y)$ or observations for all $x\in\chi$ at any observed $y$ are not compulsory.

\subsection{Properties of SNE under random sampling }
Based on the definition of $\hat{f}_x(y)$ specified by SNE, we analyze its consistency and convergence rate conditioned on a general assumption on the generation of $\mathcal{D}$. We assume that simulation observations are generated randomly and iteratively by following a distribution on each stage. To specify the distribution, first denote the probability of sampling on decision $x$ at time $t$ as $p_x^{(t)}$ (i.e.  $p_x^{(t)}>0$, $\sum_{x\in\chi}p_x^{(t)}=1$). Conditioned on $x$, we require that the marginal probability density function of sampling on $y$ exists and is denoted by $f_{Y|x}^{(t)}(y)$. Notice that this conditional distribution defined on $Y^\circ$ shares the same $\sigma$-algebra of $P_Y$. The distribution of the random sampling rule can thus be specified by:
\begin{equation}
	a^{(t)}(x, y):=p_x^{(t)}f_{Y|x}^{(t)}(y),
\end{equation} 
in the sense that the probability of sampling inside a region $V\subset\chi\times Y^\circ$ at stage $t$ can be calculated by $\sum_{x\in\chi}\int_{y: (x,y)\in V}a^{(t)}(x,y)dy$. Regardless of whether $a^{(t)}(x,y)$ on stage $t$ is derived based on last stage data $\mathcal{D}^{(t-1)}$ or predetermined for each stage $t$, we assume that each $(x,y)$ has a minimum probability to be sampled:
\begin{assumption}\label{assumption:continuous pdf}
	For a dynamic sequential random sampling rule represented by the sequence of $\{a^{(t)}(x,y)\}_t$, $$a^{(t)}(x,y)>\lambda_L,$$ $\forall t, x\in\chi, y\in Y^\circ$, for some constant $\lambda_L>0$. 
\end{assumption}

The random sampling is {\bf dynamic} if on each stage $t$, $a^{(t)}(x,y)$ is different. Otherwise, it is {\bf static}. We first show the general properties of the dynamic random sampling rule. Then we exploit special properties for optimal static random sampling and utilize them to derive an efficient dynamic random sampling rule.

\subsubsection{Analysis of SNE under Dynamic Random Sampling}
We show that for a dynamic sampling rule that satisfies Assumption \ref{assumption:continuous pdf}, $\hat{f}_{x}^{(t)}(y)\to f_x(y)$ almost surely at arbitrary point $(x,y)$ where $x\in\chi, y\in Y^\circ$, as the number of observations increases to infinity. The almost surely convergence of $\hat{f}_{x}^{(t)}(y)$ on $Y$ pointwisely for each $x$ leads to $\varpi^{(t)}$ converging to $0$. 

To prove consistency of $\hat{f}_{x}^{(t)}(y)$ at any $(x,y)$, we categorize two unwanted neighborhood situations. The first situation is when fewer samples are accumulated in the neighborhood such that $\Psi_{x, y}^{(t)}$ is smaller than expected. The second situation is when $\Psi_{x, y}^{(t)}$ is large enough yet the estimation error is too large.


To define having ``sufficient" samples in the neighborhood, we use a bound denoted as $\psi^{(t)}$ on stage $t$. Specifically, if $\Psi^{(t)}_{x, y}<\eta \psi^{(t)}$ for a specified ancillary coefficient $\eta\in(0,1)$, we regard that the neighborhood of $(x, y)$ has insufficient samples at $t$.  

\begin{lemma}
A random sampling rule satisfying Assumption \ref{assumption:continuous pdf} can guarantee that $\forall x\in\chi, y\in Y^\circ$, when $t\to\infty$,
\begin{equation}
	P\{\liminf_{t\to\infty}	\Psi_{x, y}^{(t)}\geq\lfloor\eta\psi^{(t)}\rfloor\}=1,
\end{equation}
where  $\eta$ is a fixed number in $(0,1)$ and:
\begin{equation}
	\psi^{(t)}:=(C\lambda_L)^\frac{1}{1+\xi}t^\frac{1}{1+\xi},
\end{equation}
where 
\begin{equation}
	C:=\frac{\pi^\frac{d}{2}}{\Gamma(\frac{d}{2}+1)}c^d.
\end{equation}where $\Gamma$ is the Euler's gamma function, $d$ is the dimension of $Y^\circ$ and $c$ is the constant mentioned in (\ref{definition: decreasing sequence}).
\label{Lemma:mt}
\end{lemma}
\paragraph{Remarks:} The rationale behind the construct of $\psi^{(t)}$ is to guarantee that $\psi^{(t)}$ is as large as possible while still maintaining $\Psi_{x, y}^{(t)}\geq\psi^{(t)}$ being an almost sure event when $t\to\infty$. The lemma shows that there exists a certain construct of $\psi^{(t)}$ that satisfies this. It suggests that, by definition of $\psi^{(t)}$, having enough samples in its neighborhood is guaranteed almost surely for any $(x, y)$. The need for auxiliary coefficient $\eta$ is to balance the probability of insufficient samples and significant estimation error with sufficient samples. $\eta$ will later be optimized to achieve the tightest bound possible.

On top of the almost surely guarantee of having a sufficient number of observations at any neighborhood of $(x,y)$, the estimation error caused by ``noise" also vanishes:
\begin{lemma}
If Assumption \ref{assumption: lipschitz}$\sim$\ref{assumption: independent} are satisfied , a sequential random sampling satisfying Assumption \ref{assumption:continuous pdf} can guarantee that $\forall x\in\chi, y\in Y^\circ $, when $ t\to\infty$ 
\begin{equation}
	\bar{\epsilon} _{x, y}^{(t)}\to 0\quad a.s., 
\end{equation}
where  \begin{equation}
	\bar{\epsilon} _{x, y}^{(t)}:=\frac{\sum_{i:x_{[i]}=x, y_{[i]}\in {b _{x, y}^{(t)}}}\epsilon_{[i]}}{\Psi_{x, y}^{(t)}}
\end{equation} 	\label{Lemma:error}
\end{lemma}
Combine with the proof that the error caused by ``bias" is also negligible asymptotically, we achieve the consistency proof of SNE under a general dynamic random sampling:

\begin{theorem}
\label{theorem: 1st 2nd}If Assumption \ref{assumption: lipschitz}$\sim$\ref{assumption: independent} are satisfied , a sequential random sampling satisfying Assumption \ref{assumption:continuous pdf} can guarantee that $\forall x\in\chi, y\in Y^{\circ} $, when $t\to\infty$
\begin{equation}
	{\hat{f}_{x}^{(t)}(y)}\to f_{x}(y)\quad a.s.\label{equation:first order}
\end{equation}
and
\begin{equation}
	{\hat{\sigma}_{x}^{(t)}(y)}\to \sigma_{x}(y)\quad a.s.\label{equation:second order}
\end{equation}
\end{theorem}
Based on Theorem \ref{theorem: 1st 2nd}, it is easy to see that the total decision loss by implementing $\hat{x}^{(t)}(y)$ at stage $t$ (i.e. $\varpi^{(t)}$) will decrease to $0$ eventually.
\begin{proposition}\label{proposition: consistency}
If Assumption   \ref{assumption: lipschitz}$\sim$\ref{assumption: independent} are satisfied ,a sequential random sampling satisfying Assumption \ref{assumption:continuous pdf}  can guarantee that when $t\to\infty$
\begin{equation}
	\varpi^{(t)}\to 0
\end{equation}
\end{proposition}
On top of the consistency guarantee of the decision loss, its convergence rate can also be bounded. The following theorem provides the convergence rate of conditional decision loss at each $y$, which is the backbone of our analysis of the convergence rate of the total decision loss.

\begin{theorem}\label{theorem: CDL convergence rate}
If Assumption \ref{assumption: lipschitz}$\sim$\ref{assumption: independent} are satisfied ,a sequential random sampling satisfying Assumption \ref{assumption:continuous pdf}  can guarantee that $\forall y\in Y^\circ$ such that $x^*(y)$ is unique, $\upsilon^{(t)}(y)$ converges to $0$ at least $\exp[-c_1 t^{\frac{1}{1+\xi}}]$ fast, for some $c_1\geq H_y>0$. Specifically,
\begin{equation}\label{equation: tau bound}
	\liminf_{t\to\infty}-\frac{1}{t^{\frac{1}{1+\xi}}}\ln \upsilon^{(t)}(y)\geq H_y>0,
\end{equation}
where
\begin{equation}
	H_y:=(C \lambda_L)^\frac{1}{1+\xi}\zeta(y) \left(-{W_{-1}\left(-\exp \left(-\zeta(y)-1\right)\right)}\right)^{-{\frac{1}{1+\xi}}}, \label{equation: H_y}
\end{equation}
where
\begin{equation}
	\zeta(y):= \min_{x\neq x^*(y)}\frac{\left(f_x(y)-f_{x^*}(y)\right)^2}{2 \left(\sigma _{x^*}^2(y)+\sigma _x^2(y)\right)},\label{definition: zeta}
\end{equation}
and $W_{-1}$ being the lower branch of Lambert $W$ function (i.e. $w=W_{-1}(z)$ if $we^w=z$ for $z\in[-e^{-1},0)]$).
\end{theorem}
\paragraph{Remarks:} Firstly, Theorem \ref{theorem: CDL convergence rate} suggests that the convergence rate of $\upsilon^{(t)}(y)$ is at least $O(\exp(-c_1t^{\frac{1}{1+\xi}}))$ fast if $H_y$ is well defined. Compared with the convergence rate $O(\exp(-c_1t))$ of R\&S with covariates where $Y$ is discrete (\cite{gao2019selecting},\cite{jin2019optimal}), the loss of convergence rate is caused by the shrink neighborhood rule such that the neighborhood sample size only increases $t^\frac{1}{1+\xi}$ fast exactingly. Although $\xi$ can be set as an arbitrarily small number to make $t^{\frac{1}{1+\xi}}$ as close as $t$, the $\xi$ cannot be $0$ while maintaining consistency of the estimator. Secondly, some remarks on $H_y$. The constant $C:=\frac{\pi^\frac{d}{2}}{\Gamma(\frac{d}{2}+1)}c^d$ is a function w.r.t. to the dimension $d$ and constant $c$ mentioned in \ref{definition: decreasing sequence} . In principle (e.g. $d>5$), the larger the $d$, the smaller the $C$. To keep the convergence rate, it is advised to set $c$ larger corresponding to a larger $d$. Since $c$ can always be adjusted according to $d$, the convergence rate of decision loss does not suffer from a high dimension of $Y$. The global lower bound for $a^{(t)}(x,y)$ (i.e. $\lambda_L$) also contributes to the convergence rate. The larger the $\lambda_L$, the faster the convergence. The last critical element in $H_y$ is $\zeta(y)$. $\zeta(y)$ represents the level of difficulty to differentiate non-optimal decisions (e.g. $x\neq x^*(y)$) from $x^*(y)$. Notice that the impact of $\zeta(y)$ to $H_y$ is monotone increasing. This means the harder it is to differentiate the optimal decision, the slower the convergence rate will be. It is also worth noticing that for context that has multiple conditional optimal $x$, $\zeta(y)$ is not well defined. We call them Multi-Optimum (MO) states. The existence of MO prohibits a convergence rate bound for $\varpi^{(t)}$. However, a bound can be derived on a compact subset of $Y^\circ$ if MO is excluded:  
\begin{proposition}\label{proposition: subset rate}

If Assumption  \ref{assumption: lipschitz}$\sim$\ref{assumption: independent} are satisfied, a sequential random sampling satisfying Assumption \ref{assumption:continuous pdf} can guarantee that for any compact subset $Y_{sub}\subset Y^\circ$ such that $Y_{sub}\cap MO=\emptyset$, $\varpi^{(t)}$ defined on this subset:
$$\varpi_{Y_{sub}}^{(t)}:=\int_{y\in Y_{sub}}\upsilon^{(t)}(y)d\mu_Y(y)$$ converges to $0$ at least $\exp[-c_2t^{\frac{1}{1+\xi}}]$ fast, for some $c_2\geq H_{Y_{sub}}>0$. Specifically: 
\begin{equation}
	\liminf_{t\to\infty}-\frac{1}{t^{\frac{1}{1+\xi}}}\ln \varpi_{Y_{sub}}^{(t)}(y)\geq H_{Y_{sub}}>0,
\end{equation}
where
\begin{equation}
	H_{Y_{sub}}:=\inf_{y\in Y_{sub}}H_y,\label{equation: error on subset}
\end{equation}
and $H_y$ defined in (\ref{equation: H_y})
\end{proposition}

\paragraph{Remarks:} From proposition \ref{proposition: subset rate}, we conclude that the convergence rate of $\varpi^{(t)}$ defined on any compact subset $Y_{sub}$ such that $Y_{sub}\cap MO=\emptyset$ is lower bounded by the slowest convergence rate of $\upsilon^{(t)}(y)$ on the set. This bound is invariant to the likelihood distribution of $y$ (i.e. $P_Y$). This invariance to the likelihood of the contexts makes the bound very stable and robust. Notice that according to this bound, the convergence speed is linked to $\lambda_L$. At a glance, it seems to be optimal to allocate the budget uniformly across $\chi\times Y_{sub}$. This will lead to the static uniform distribution being optimal. However, we shall analyze the static random sampling rule and show that there exists optimal static random sampling which is not uniform sampling.

\subsubsection{Analysis of SNE under Static Random Sampling}
The analysis of the general dynamic random sampling suggests that a uniform sampling may be optimal. In this section, we show that if the sampling rule is static, a tighter bound for the rate function of $\varpi_{Y_{sub}}^{(t)}$ can be achieved. In terms of this tighter bound, there exists an optimal static random sampling rule. This inspires us to derive a dynamic rule that can converge to this optimal static rule.

For a static rule, a unique $a(x,y)$ function is applied on each stage $t$. To achieve convergence, we assume that $\forall x\in\chi$, $a(x,y)$ is strictly positive and continuous on $Y^\circ$: 
\begin{assumptionbis}{assumption:continuous pdf}\label{assumption:continuous pdf static}
For a static random sampling rule represented by $a(x,y)$, $\forall x\in\chi, y\in Y^\circ$, $$a(x,y)>0.$$
In addition, $\forall  x\in\chi$, $a(x,y)$ is continuous on $y\in Y^\circ$. 
\end{assumptionbis}

\begin{lemma}
A static random sampling rule satisfying Assumption \ref{assumption:continuous pdf static} can guarantee that $\forall x\in\chi, y\in Y^\circ$, when $t\to\infty$,
\begin{equation}
	P\{\liminf_{t\to\infty}	\{\Psi_{x, y}^{(t)}\geq\lfloor\eta_{x,y}\psi_{x,y}^{(t)}\rfloor\}\}=1,
\end{equation}
where  $\eta_{x,y}$ is a fixed number in $(0,1)$ associated with each $(x,y)$ and:
\begin{equation}
	\psi_{x, y}^{(t)}:=\max\{m\in[1,t],m\in\mathbb{N}:[Ca_{\inf,x,y,m}]^\frac{1}{1+\xi}t^\frac{1}{1+\xi}\geq m\},
\end{equation}
where 
\begin{equation}
	\begin{split}		
		C:=&\frac{\pi^\frac{d}{2}}{\Gamma(\frac{d}{2}+1)}c^d\\
		a_{\inf,x,y,m}:=&\inf_{y'\in B(r_m;y)}a(x,y').
	\end{split}
\end{equation}where $\Gamma$ is the Euler's gamma function, $d$ is the dimension of $Y^\circ$ and $c$ is the constant mentioned in (\ref{definition: decreasing sequence}).
\label{Lemma:mt static}
\end{lemma}
\paragraph{Remarks:} Comparing to the analysis on dynamic random sampling rule, a $\psi_{x,y}^{(t)}$ can now be associated with each $x,y$. The existence of $\psi_{x, y}^{(t)}$ for any $t$ and $(x,y)$ can be secured by setting $c$ large enough so that $[Ca_{\inf,x,y,2}]^\frac{1}{1+\xi}\geq 2$. This will guarantee that after initializing the static random sampling rule by observing $2$ contexts for each $x\in\chi$, at least $b_2^{(2)}(x,y)$ is activated for any $(x,y)$. 

\begin{theorem}\label{theorem: CDL convergence rate static}
If Assumption \ref{assumption: lipschitz}$\sim$\ref{assumption: independent} are satisfied, a static random sampling specified by $a(x,y)$ satisfying Assumption \ref{assumption:continuous pdf static} can guarantee that $\forall y\in Y^\circ$ such that $x^*(y)$ is unique, $\upsilon^{(t)}(y)$ converges to $0$ at least $\exp[-c_3 t^{\frac{1}{1+\xi}}]$ fast for some constant $c_3\geq H'_y$. Specifically, given $y\in Y^\circ$
\begin{equation}
	\begin{split}
		\lim_{t\to\infty}-\frac{1}{t^{\frac{1}{1+\xi}}}\ln \upsilon^{(t)}(y)\geq & H'_y(\{a(x,y)\}_{x\in\chi})\\
		:=&C^\frac{1}{1+\xi}\min_{x\in\chi}  \left\{ \zeta_{x}(y) \left(-W_{-1}\left(-\exp[{-    \zeta_{x}(y)   -1}]\right)\right)^{-\frac{1}{1+\xi}} [a(x,y)]^\frac{1}{1+\xi}\right\} ,   \\
	\end{split}\label{equation: H_y static}
\end{equation}
where
\begin{equation}
	\zeta_{x}(y):=\begin{cases}
		\frac{(f_x(y)-f_{x^*}(y))^2}{4 \sigma_{x}^2(y)}& x\neq x^*(y)\\
		\frac{(   \min_{x\neq x^*(y)}  f_x(y)-f_{x^*}(y))^2}{4\sigma_{x^*}^2(y)} & x=x^*(y)
	\end{cases},
\end{equation}	
and $W_{-1}$ being the lower branch of Lambert $W$ function (i.e. $w=W_{-1}(z)$ if $we^w=z$ for $z\in[-e^{-1},0)]$).
\end{theorem}
\paragraph{Remarks:} Comparing to $H_y$ in Theorem \ref{theorem: CDL convergence rate}, $H'_y$ is associated with $a(x,y)$ rather than its global lower bound $\lambda_L$. This shows that $H'_y$ is a tighter bound than $H_y$.

The rate bound for $\upsilon^{(t)}(y)$ leads to the rate bound of $\varpi^{(t)}$ on compact set $Y_{sub}\subset Y^\circ$ excluding MO contexts. This is shown in the following Proposition.

\begin{proposition}\label{proposition: subset rate static}
If Assumption \ref{assumption: lipschitz}$\sim$\ref{assumption: independent} are satisfied , a static random sampling specified by $a(x,y)$ satisfying Assumption \ref{assumption:continuous pdf static} can guarantee that for any compact subset $Y_{sub}\subset Y^\circ$ such that $Y_{sub}\cap MO=\emptyset$, $\varpi^{(t)}$ defined on this subset:
$$\varpi_{Y_{sub}}^{(t)}:=\int_{y\in Y_{sub}}\upsilon^{(t)}(y)d\mu_Y(y)$$ converges to $0$ at least $\exp[-c_1t^{\frac{1}{1+\xi}}]$ fast, for some $c_1\geq H'_{Y_{sub}}>0$. Specifically: 
\begin{equation}
	\liminf_{t\to\infty}-\frac{1}{t^{\frac{1}{1+\xi}}}\ln \varpi_{Y_{sub}}^{(t)}(y)\geq H'_{Y_{sub}}(\{a(x,y)\}_{x\in\chi, y\in Y_{sub}})>0,
\end{equation}
where
\begin{equation}
	H'_{Y_{sub}}(\{a(x,y)\}_{x\in\chi, y\in Y_{sub}}):=\inf_{y\in Y_{sub}}{H'}_y(\{a(x,y)\}_{x\in\chi}),\label{equation: error on subset static}
\end{equation}
with $H'_y$ defined in (\ref{equation: H_y static})
\end{proposition}

\paragraph{Remarks:}Compared to Proposition \ref{proposition: subset rate}, the bound shown in Proposition \ref{proposition: subset rate static} is related to $a(x,y)$. This leads us to the idea of maximizing $H'_{Y_{sub}}$ on a large enough $Y_{sub}$, which is equivalent to the following optimization problem:
\begin{equation}
\label{equation: optimal budget allocation}
\begin{split}
	\max_{a(x,y),\forall x\in \chi, y\in Y^\circ} &\  H'_{Y_{sub}}\\
	s.t. &\sum_{x\in\chi}\int_{Y_{sub}}a(x, y)dy=1.
\end{split} 
\end{equation}
Although the above optimization problem relies on full information of $f_x(y)$ and $\sigma_x(y)$ on $\chi\times Y^\circ$, by satisfying Assumption \ref{assumption:continuous pdf}, Theorem \ref{theorem: 1st 2nd} shows that a dynamic rule can guarantee that $\hat{f}_x(y)$ and $\hat{\sigma}_x(y)$ from SNE is consistent. By substituting $f_x(y)$ and $\sigma_x(y)$ by $\hat{f}_x(y)$ and $\hat{\sigma}_x(y)$ and resolving the optimization iteratively, we could converge to the optimal static allocation.

\section{Contextual Optimizer through Neighborhood Estimation}
\label{section:methodology_Allocation}
In this section, we illustrate our allocation rule called Contextual Optimizer through Neighborhood Estimation(CONE). The idea is a dynamic random sampling rule that satisfies Assumption \ref{assumption:continuous pdf} while converging to the static optimal distribution derived from (\ref{equation: optimal budget allocation}).
\subsection{Optimal static budget distribution and its approximation}
We have the following theorem stating the optimality of $a(x,y)$ in terms of $H'_{Y_{sub}}$:
\begin{theorem}
\label{Theorem: optimal budget allocation}
The following is an optimal solution of problem (\ref{equation: optimal budget allocation}):
\begin{equation}
	a^*(x,y)=\frac{\lambda \left(-W_{-1}\left(-\exp[{-    \zeta_{x}(y)   -1}]\right)\right)}{C[\zeta_{x}(y)]^{1+\xi} }\quad\forall x\in\chi, y\in Y_{sub}\label{equation: optimal allocation}
\end{equation}
where
\begin{equation}
	\begin{split}
		\zeta_{x}(y):=&\begin{cases}
			\frac{(f_x(y)-f_{x^*}(y))^2}{4 \sigma_{x}^2(y)}& x\neq x^*(y)\\
			\frac{(   \min_{x\neq x^*(y)}  f_x(y)-f_{x^*}(y))^2}{4\sigma_{x^*}^2(y)} & x=x^*(y)
		\end{cases},\\
		\lambda=&\frac{1}{\sum_{x\in\chi}\int_{Y_{sub}}	\frac{\left(-W_{-1}\left(-\exp[{-    \zeta_{x}(y)   -1}]\right)\right)}{C [\zeta_{x}(y)]^{1+\xi} } dy},
	\end{split}
\end{equation}
and $W_{-1}$ being the lower branch of Lambert $W$ function (i.e. $w=W_{-1}(z)$ if $we^w=z$ for $z\in[-e^{-1},0)]$).
\end{theorem}
\paragraph{Remarks:}Notice that $a(x,y)$ is monotone decreasing w.r.t. $\zeta_{x}(y)$.  From the formulation, $\zeta_{x}(y)$ is the ratio between the expected performance difference to the optimal and the heteroscedastic noise. It represents how easy it is for this $(x,y)$ to be ranked correctly conditioned on $y$. The easier it is, the less budget should be spent on its neighborhood. A good feature of this solution is that although it is derived on a subset of $Y^\circ$ with no MO contexts, the rule itself can be implemented on $Y$ as long as the $\lambda$ is well defined. If this is the case, letting $a(x,y)=\infty$ at some points won't cause an issue. However, empirically, a $Y_{sub}$ that covers regions too close to MO contexts will lead to a very large $a(x,y)$. This could result in excessively greedy behavior, anticipated to be detrimental during the early stages when only a limited number of observations have been accumulated. Unstable estimations could potentially exhaust the budget on regions initially estimated to be proximal to the multi-objective (MO) but are, in fact, easily distinguishable upon further observation. To rectify this problem, a global upper bound will be imposed to truncate $a(x,y)$, limiting any value exceeding the bound.

Our algorithm, Contextual Optimizer through Neighborhood Estimation (CONE), allocates observations based on distribution $a^*(x,y)$ formulated in (\ref{equation: optimal allocation}). Since the real $f_x(y)$ and $\sigma_{x}(y)$ are unknown, on each stage, they are estimated by $\hat{f}^{(t)}$ and $\hat{\sigma}^{(t)}$ by SNE. Thus on stage $t$:
\begin{equation}
a^{(t)}(x,y)=\frac{\lambda \left(-W_{-1}\left(-\exp[{-    \hat{\zeta}_{x}^{(t)}(y)   -1}]\right)\right)}{C[\hat{\zeta}_{x}^{(t)}(y)]^{1+\xi} }\quad\forall x\in\chi, y\in Y_{sub}
\end{equation}
where
\begin{equation}
\begin{split}
	\hat{\zeta}_{x}(y):=&\begin{cases}
		\frac{(\hat{f}_x^{(t)}(y)-\hat{f}_{\hat{x}^{(t)}}^{(t)}(y))^2}{4 [\hat{\sigma}_{x}^{(t)}(y)]^2}& x\neq \hat{x}(y)\\
		\frac{(   \min_{x\neq \hat{x}(y)} \hat{f}_x^{(t)}(y)-\hat{f}_{\hat{x}^{(t)}}^{(t)}(y))^2}{4[\hat{\sigma}^{(t)}_{\hat{x}^{(t)}}(y)]^2} & x=x^*(y)
	\end{cases},
\end{split}
\end{equation} Since only the relative ratio is important in empirical sampling, $\lambda$ and $C$ can be ignored when applying Rejection Sampling or Markov Chain Monte Carlo methods. We can thus focus on 
$$\beta^{(t)}(x,y):=\frac{ \left(-W_{-1}\left(-\exp[{-    \hat{\zeta}_{x}^{(t)}(y)   -1}]\right)\right)}{[\hat{\zeta}_{x}^{(t)}(y)]^{1+\xi} }.$$
To satisfy Assumption \ref{assumption:continuous pdf} so that CONE can guarantee convergence, a lower bound $\lambda_L$ for $\beta^{(t)}(x,y)$ is applied on all stages. To prevent greediness near the MO context, we impose a global upper bound $\lambda_U$ on $\beta^{(t)}(x,y)$. So any $\beta^{(t)}(x,y)$ smaller than $\lambda_L$ will be truncated up to $\lambda_L$ and any that is larger than $\lambda_U$ will be truncated down to $\lambda_U$. 

In summary, CONE allocates budget based on estimated optimal static distribution. By satisfying Assumption \ref{assumption:continuous pdf}, CONE guarantees to provide the proxy with decision loss converges to $0$ at near-to-optimal rate. For a version of sequential budget allocation that uses Rejection Sampling to realize $a^{(t)}(x,y)$, one can refer to the supplement.

\section{Numerical Examples}
\label{section: numericals}
In this section, we present empirical evidence underscoring our two principal contributions: the development of a consistent surrogate, the Shrinking Neighborhood Estimation (SNE), which assures a rapid convergence rate on estimation error, and the introduction of Rate-Efficient-Sampling (RES), a methodology constructed upon the SNE framework, devised to further expedite this rate of convergence. The efficacy of the RES and SNE combination is substantiated via two case studies. The first case is a challenging one-dimensional toy example, and the second case is a real-world staffing scenario in a call center, where the complex interdependencies of the queuing system necessitate reliance on simulation modeling.

To differentiate the contribution of SNE and RES to performance, we compare three algorithms:
\begin{enumerate}
\item USKrig: This algorithm distributes the budget uniformly on $\chi\times Y$ space. When the budget is depleted, the normal Kriging model is applied to data from the same $x$ on $Y$ space to estimate performance on unobserved $y$ point. The conditional optimal decision is identified by the ranking of the estimation given by Kriging at $y$. To enable using Kriging in heteroscedastic noise, the budget is allocated in batches with a fixed size for each new point and the batch variance is applied to estimate the variance of the local noise. When the budget is depleted, the Kriging model with estimated variance at each observed point added to the diagonal of the covariance matrix is applied. In this research, we set the batch size as $5$.
\item USSNE: This algorithm distributes the budget uniformly on $\chi\times Y$ space. When the budget is depleted, Shrinking Neighborhood Estimation is applied to data from the same $x$ on $Y$ space to estimate performance on unobserved $y$ point. The conditional optimal decision is identified by the ranking of the estimation given by SNE at each $y$. SNE does not require replications, and each point only gets one observation.
\item CONE: This algorithm distributes the budget optimally based on SNE. When the budget is depleted, Shrinking Neighborhood Estimation is applied to data from the same $x$ on $Y$ space to estimate performance on unobserved $y$ point. The conditional optimal decision is identified by the ranking of the estimation given by SNE at each $y$. SNE does not require replications, and each point only gets one observation.
\end{enumerate}
On one hand, the comparison between USKrig and USSNE can demonstrate the capability of SNE as a superior estimator in a heteroscedastic scenario. On the other hand, the comparison between USSNE and CONE can further show the superior performance of Rate Efficient Sampling as compared to Uniform Sampling.
\subsection{Toy case}
\begin{figure}[!h]
\centering
\includegraphics[scale=0.4]{ 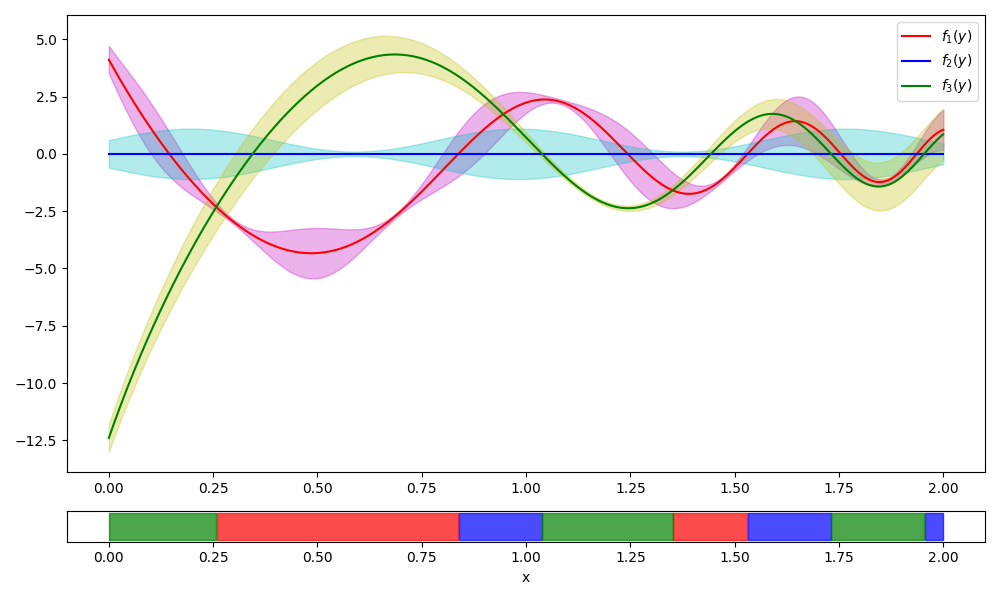}
\caption{The plot of $f_{x}(y)$ in 1-d benchmark problem with the red, blue, green indicating $x=1, 2, 3$ respectively. The shade that covers each plot shows $\pm\sigma_x(y)$ region. The strip below color-codes the conditional optimal decision $x^*(y)$ under each state $y$. }\label{Fig: 1d_real_f}
\end{figure}

\begin{figure}[!h]
\centering
\includegraphics[width=0.7\textwidth]{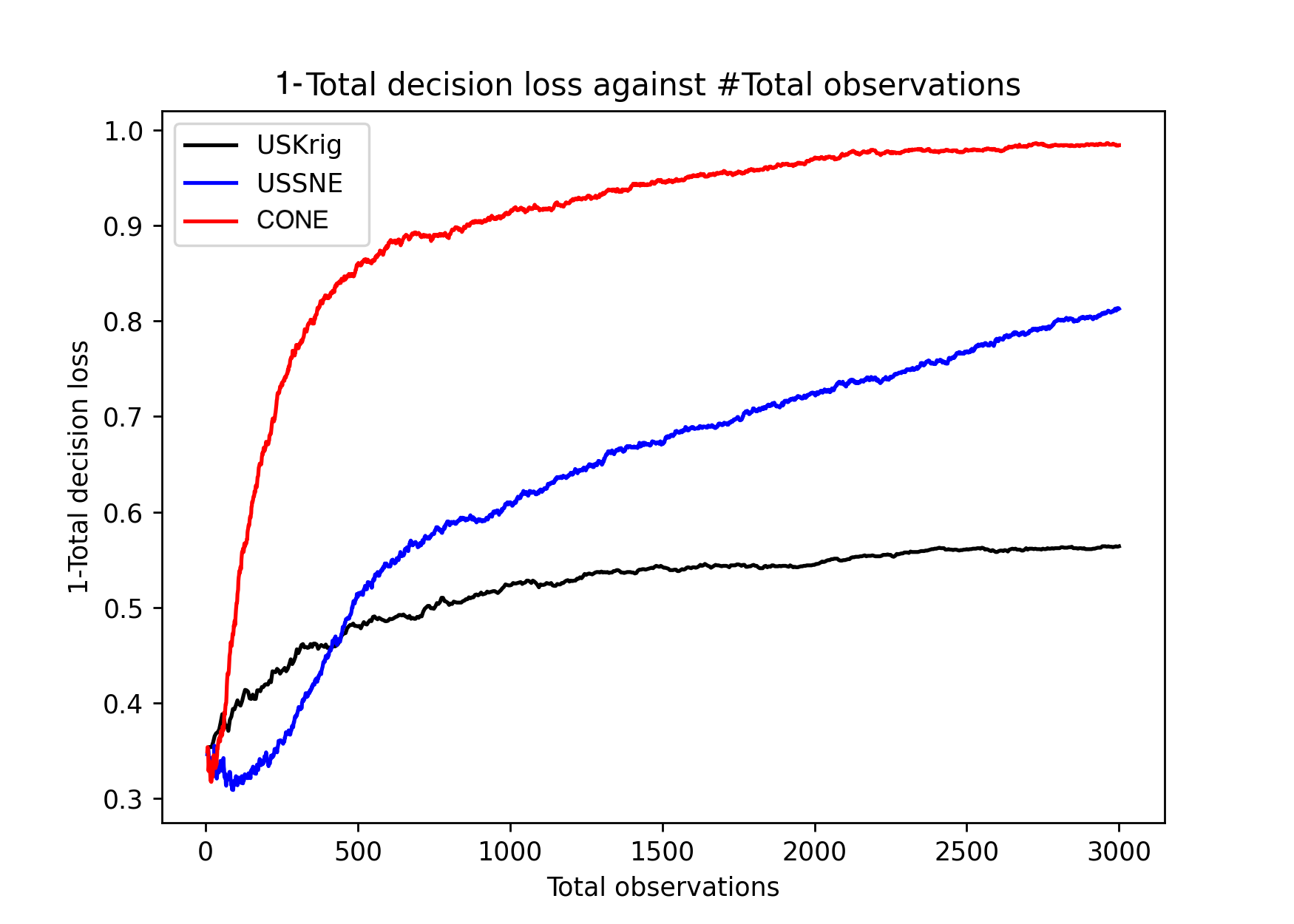}
\caption{Plots of $1-\varpi^{(t)}$ against the number of observations $t$ for three algorithms in solving 1-d benchmark problem. The higher the value, the better the algorithm. The result is the average among $100$ independent trials.}\label{Fig: 1d_TDL}
\end{figure}

We test on a 1-d $Y$ benchmark problem. It is fabricated in a fashion so that the following three features challenging the correct identification of the conditional optimality are reflected:
\begin{enumerate}
\item The conditional optimal decision changes under different $y$ and the distribution of MO contexts is uneven. For regions where the conditional optimal decisions change frequently, it is more difficult to identify the conditional optimal decision on them.
\item The performance difference between any two decisions should be different across $Y$. For those regions where the conditional best and the conditional second best perform similarly, it is more difficult to identify the conditional optimal decision on them.
\item Heteroscedastic noise. The variance $\sigma_x(y)$ is different and fluctuate on $\chi\times Y$.
\end{enumerate}
The following is the designed benchmark problem:
\begin{align}
\tilde{G}_x(y)&\sim\mathcal{N}(f_{x}(y), \sigma_x(y))\\
f_{x}(y)&=\begin{cases}
	\frac{10}{(y+1)^2}\sin(e^{y+1})&x=1	\\
	0&x=2\\
	-\frac{10}{(y+0.8)^2}\sin(e^{y+0.8})&x=3
\end{cases}\quad y\in [0, 2]\\
\sigma_x(y)&=\begin{cases}
	0.5 (\sin (16 y)+1.2)&x=1	\\
	0.5 (\sin (8 y)+1.2)&x=2\\
	0.5 (\sin (4 y)+1.2)&x=3
\end{cases} \quad y\in [0, 2]
\end{align}We assume that $\mu_Y$ is uniform distribution on $Y$. The plot of $f_{x}(y)$ is shown in Figure \ref{Fig: 1d_real_f}. 
We run $100$ independent macro trails for each of the three algorithms to solve the toy case problem given $3000$ total observation budget. For each macro trail and a given algorithm, at stage $t$, the total decision loss $\varpi^{(t)}$ is estimated by Monte Carlo sampling $20$ points on $Y$, and calculate the percentage of $y$ that satisfies $\hat{x}^{(t)}(y)\neq x^*(y)$. We then average the Monte Carlo estimated $\varpi^{(t)}$ among $100$ trails. The plot of the Monte Carlo estimation on $\varpi^{(t)}$ averaged among $100$ macro trails are shown in Figure \ref{Fig: 1d_TDL}.

As shown in the plot, after $3000$ observations accumulated, $1-\varpi^{(t)}$ of USKrig stabilized at around $0.55$. This matches our anticipation of the normal Kriging method where it struggles to differentiate whether the local fluctuation of the observed performance is caused by noise or fluctuation of the underlining functions, leading to a poor estimation and incorrect identification. On the contrary, $1-\varpi^{(t)}$ of SNE methods shows a clear trend to increase to $1$, which matches our previous proof that SNE surrogate is consistent. The comparison between USSNE and CONE shows that by recursively optimizing based on SNE estimated convergence rate of Total Decision Loss, the performance is significantly improved. With $3000$ single observation accumulated, CONE is expected to achieve almost $100\%$ correct identification.

\subsection{Call center staffing problem}
In the context of healthcare operations, hospital call centers serve crucial functions such as scheduling appointments (AP) and addressing general enquiries (GE). Efficient staffing for these functions, particularly during peak morning shift hours, is a significant operational challenge. This section presents an in-depth analysis of this problem, demonstrating the promising application of CONE method to address this issue.

The hospital call center under scrutiny operates two lines: the AP and GE. There are $60$ staff members. Each of them can serve both lines. The staffing level for these lines remains fixed during a shift, which spans a continuous four-hour period, with our focus being the morning shift from 8:00 to 12:00. A notable hurdle in this system is call abandonment, which transpires when all operators are occupied and callers disconnect after waiting beyond their patience threshold. To model the behavior of callers, the arrival rate and their maximum patience are modeled independently for each half-hour slot and each line. That accounts for $16$ slot and line combinations.  From existing data, the empirical distribution of maximum caller patience is known and fixed for each of these combinations.

\begin{figure}[!h]
\centering
\includegraphics[width=0.4\textwidth]{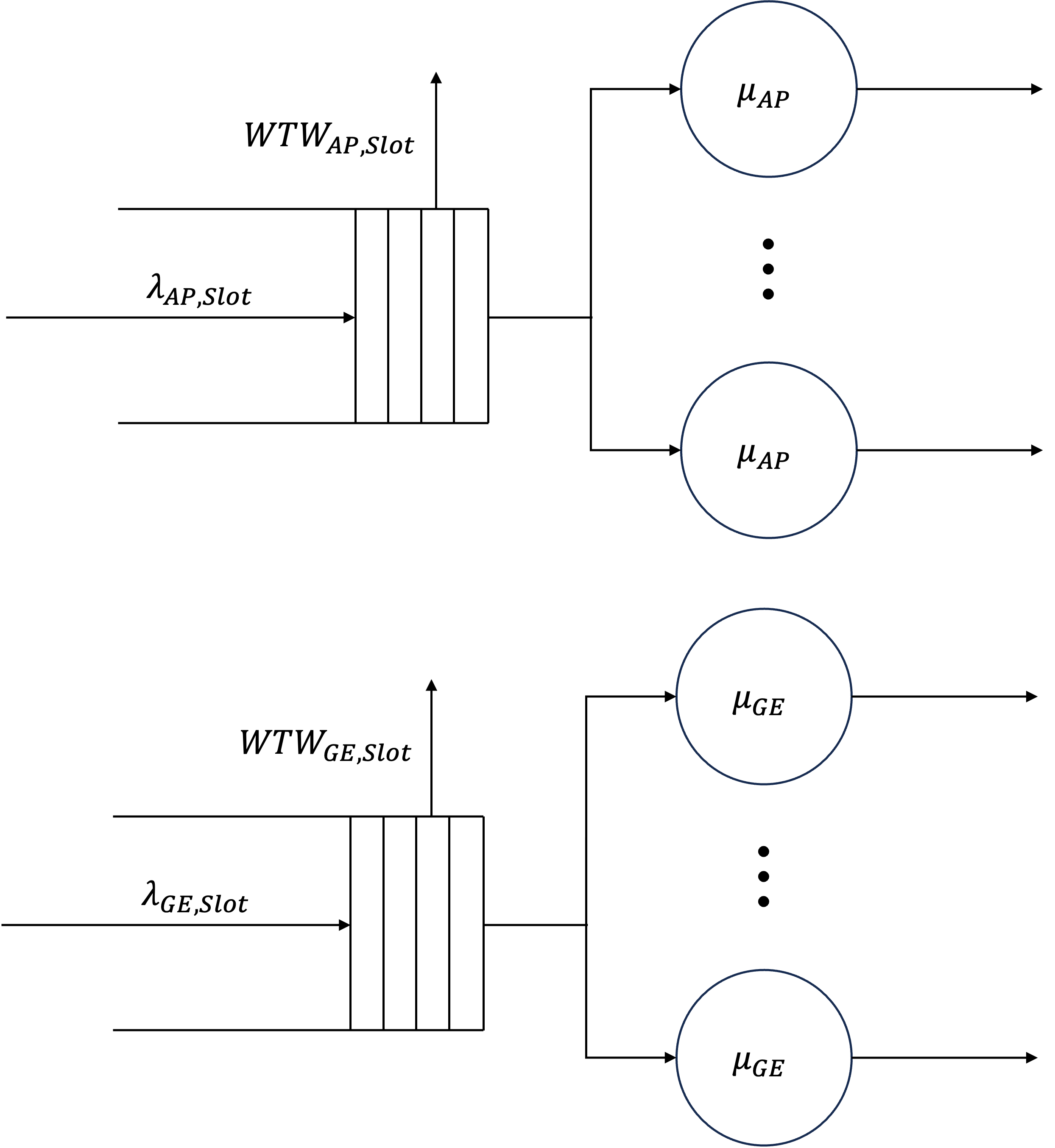}
\caption{An illustration depicts the call center's queuing system with each caller's empirically derived willingness-to-wait triggering abandonment. Both the arrival rate $\lambda$ and the Willingness-To-Wait distribution are set to adjust every 30 minutes slot.}\label{Fig: Q}
\end{figure}

Given these parameters, the aim is to determine an optimal fixed split of the $60$ staff between the AP and GE lines for the morning shift that accounts for the forecast arrival rate for each slot and line combination, with the overarching objective of minimizing the overall abandonment rate.

\begin{figure}[!h]
\centering
\includegraphics[width=0.8\textwidth]{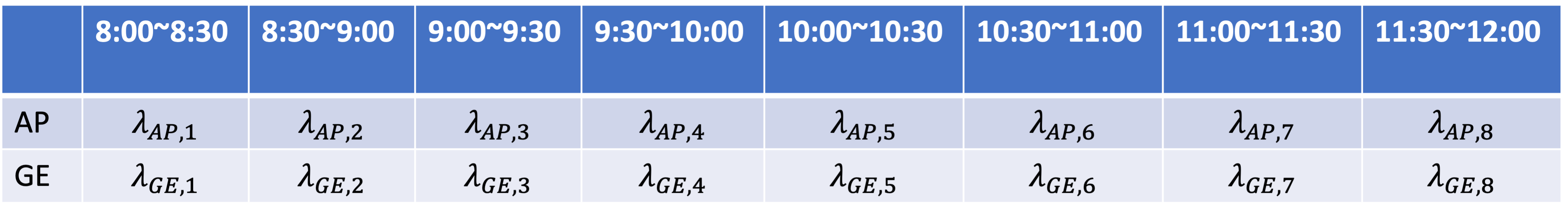}
\caption{A table showing $16$ dimension context variable.}\label{Fig: slot}
\end{figure}

\begin{figure}[!h]
\centering
\includegraphics[width=0.5\textwidth]{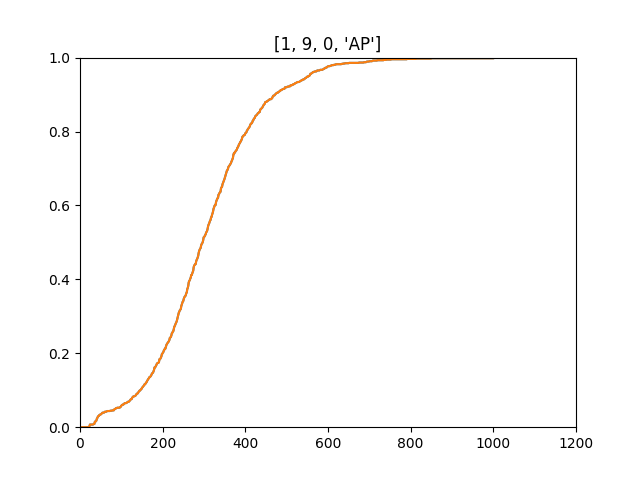}
\caption{The empirical CDF of the distribution of caller's maximum willingness-to-wait before been served on Monday morning during 9:00$\sim$9:30.}\label{Fig: WTW}
\end{figure}
The main challenge of this problem is that each half-hour arrival rate prediction is made each early morning, and the call center does not have sufficient time to identify the optimal split conditioned on each arrival profile. In addition, due to the correlation of abandonment rate between each half-hour plus the complexity of caller patience, the relation between the overall abandonment rate and a specific staffing split can be regarded as a black box. 

To rigorously link this problem with our notation. The decision $x\in\chi$ corresponds to each valid way to split $60$ staff. By holding the requirements such that each line should have at least $5$ staff, there are $51$ ways of splitting in total. The context $y\in Y\subset \mathbb{R}^{16}$ is a $16$ dimension vector specifying the arrival rate for AP and GE line for each of $8$ half-hour slots within the morning shift.  $Y$ is a rectangle, with each dimension constrained by the maximum and minimum arrival rate seen in history. We assume that the likelihood of $y$ is uniform on $Y$. The performance $f_x(y)$ is the expected overall abandonment rate in the morning shift conditioned on making a specific split under arrival profile $y$. The goal is to train a proxy $\hat{x}(y)$ that approximates $x^*(y)$ which indicates the optimal split minimizing the abandonment rate conditioned on each arrival profile.

To solve this problem, a simulation is established to model the dynamics of the queuing system based on arrival and caller patience profile. Since $f_x(y)$ is unknown ahead when estimating $\varpi$, we use average performance among $10$ additional simulation run of reach $x\in\chi$ conditioned on $y$ to estimate the real best. We implement USKrig, USSNE, and CONE and compare their $1-\varpi$ after $3000$ simulation runs. We repeat this whole process for $10$ times and show the mean and variance of $1-\varpi$.

\begin{figure}[!h]
\centering
\includegraphics[width=0.7\textwidth]{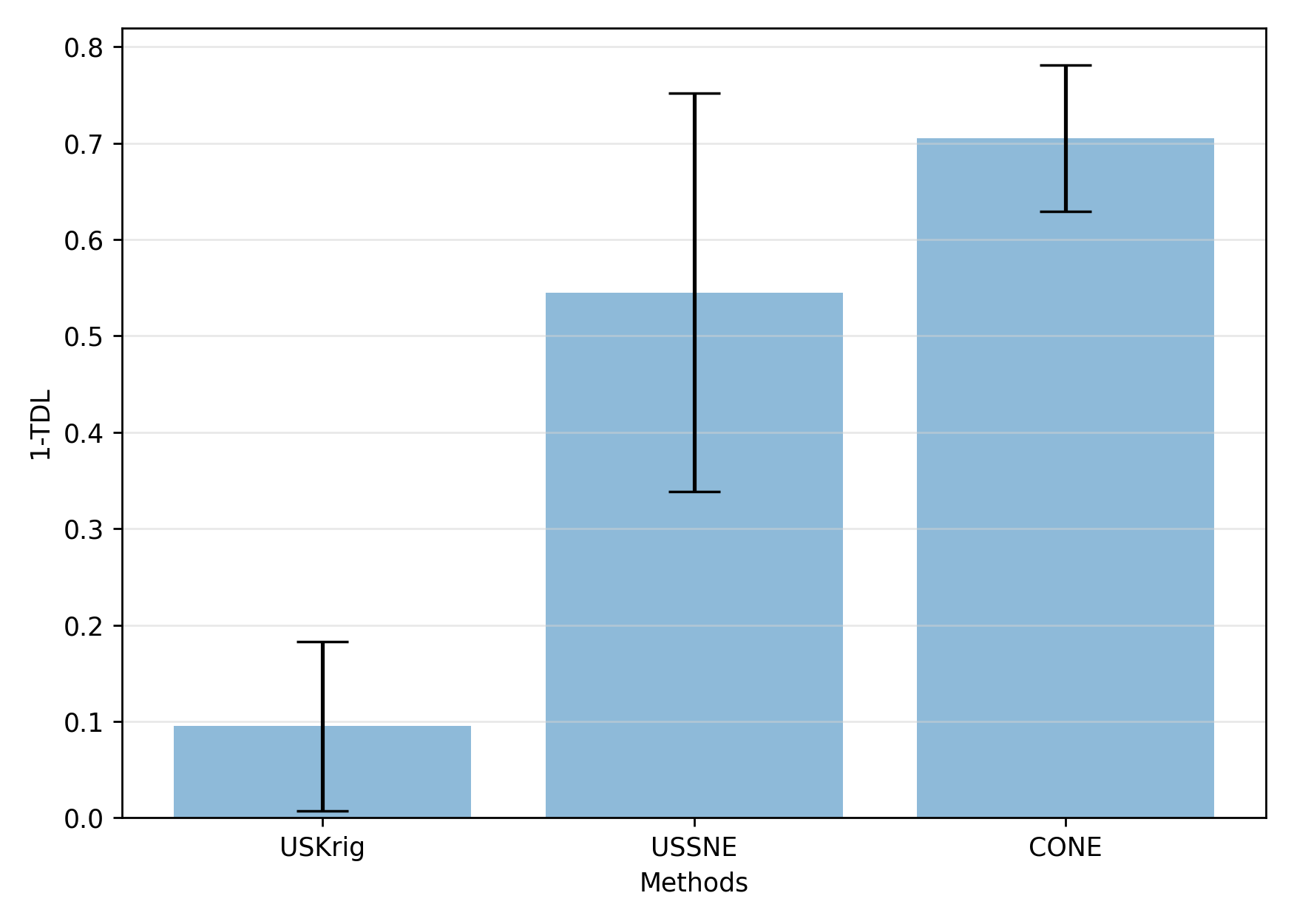}
\caption{Monte Carlo estimated $1-\varpi^{(t)}$ for three algorithms after $3000$ simulation runs. The higher the value, the better the algorithm. The range on the bar shows one standard deviation of $1-\varpi^{(t)}$ across $10$ independent meta-experiments.}\label{Fig: Usecase}
\end{figure}

As shown in Figure \ref{Fig: Usecase}, the empirical evidence demonstrates that $\hat{x}(y)$, trained via CONE, delivers the most optimal performance, closely followed by USSNE. Notably, as both USKrig and USSNE predominantly employ Uniform Sampling, the comparison between them validates that a mere surrogate model shift from Kriging to SNE results in a considerable enhancement in successful identification: a rise from a modest 10\% to an impressive 55\%. Building upon this significant improvement, the incorporation of Rate-Efficient Sampling fortifies the utility of each simulation run, endowing CONE with an additional 15\% advantage over USSNE. This outcome affirms the efficacy of the alliance between RES and SNE in resolving high-dimensional heteroscedastic problems. 

\section{Conclusion}
In the face of high-dimensional, heteroscedastic problems encountered within the scope of operations research, this paper delivers efficient prescriptive solutions that address the identified gaps in the literature.

Our first contribution lies in the proposition of the Shrinking Neighborhood Estimation method. The uniqueness of this method rests in its adaptive capability to heteroscedastic noise, eradicating the necessity of compulsory replications at each point $(x, y)$. Importantly, through this approach, consistent estimates on both the expected performance and variance are realized, ultimately retrieving the real conditional optimal decision for an arbitrary state $y$. The theoretical underpinning, coupled with the demonstrated results, accentuates the method's efficacy and its potential to be generalized across similar problem domains.

Secondly, we proposed an innovative sequential rule based on an optimal convergence rate. This method addresses a critical aspect of budget allocation and steers the budget toward those $(x, y)$ combinations that significantly impact the overall decision loss. Our comparative study with benchmark algorithms elucidates the superior performance of this methodology in terms of efficiency and results, contributing an effective tool for budget allocation in similar contexts.

Lastly, the paper manifests a real-world application of the proposed methods by addressing a staffing problem in a hospital call center. Here, the ultimate goal was to provide an optimal staffing allocation considering caller arrival and patient profile. The results suggest promising potential for the proposed methods in facilitating high-quality prescriptive analysis for problems of similar nature, thereby offering direct practical implications.

In conclusion, our study presents substantial contributions to the field of operations research, delivering both methodological advancements and practical insights. By successfully addressing the complexities and challenges associated with high-dimensional, heteroscedastic problems, this work opens up avenues for future research and potential improvements in the operational efficiency of various industries.
\label{section: conclusion}

\ACKNOWLEDGMENT{The authors gratefully acknowledge the support by Centre for Next Generation Logistics (C4NGL), Centre of Excellence in Modelling and Simulation for Next Generation Ports (C4NGP) and National University Health System (NUHS).}


\bibliographystyle{informs2014} 
\bibliography{mybib.bib} 



\newpage

\section{Proof of statements}



\subsection{Lemma \ref{Lemma:mt}}
\proof{Proof}
Since the discussion is for each $x, y$, we omit specification of $x, y$ for $\Psi^{(t)}_{x,y}$ if no confusion is caused. We prove by using Borel-Cantelli Lemma, which is to show that $\forall \eta\in(0,1)$, 
\begin{align*}
\sum_{t=1}^\infty P\{\Psi^{(t)} <\lfloor\eta\psi^{(t)} \rfloor\}<\infty.
\end{align*}
Notice that by definition:
\begin{equation}
\psi^{(t)}:=(C\lambda_L)^\frac{1}{1+\xi}t^\frac{1}{1+\xi},\label{definition: psi}
\end{equation}
where \begin{equation}
C:=\frac{\pi^\frac{d}{2}}{\Gamma(\frac{d}{2}+1)}c^d,
\end{equation}
and $\lambda_L$ is the global lower bound for $a$. On stage $t$, denote the probability such that a point $(x',y')$ randomly generated by $a^{(t)}(x,y)$ falls into the neighborhood of $(x,y)$ with radius of $r_m$ (e.g. $x'=x, y'\in B(r_m;y)$) as $p_m^{(t)}$ (again, $x,y$ is omitted in this notation).  Notice that in the following discussion, $m$ is allowed to take a continuous value. For $t\geq T^*$, we have the following lower bound for $p_{\psi^{(t)}}^{(t)}$:
\begin{equation}
p_{\eta\psi^{(t)}}^{(t)}\geq\frac{\pi^\frac{d}{2}}{\Gamma(\frac{d}{2}+1)}r_{{\eta\psi^{(t)}}}^d\lambda_L=C(\eta\psi^{(t)})^{-\xi}\lambda_L\label{inequality: p bound}
\end{equation}
The inequality holds as $p_{\eta\psi^{(t)}}^{(t)}$ is the integral of $a^{(t)}(x,y)$ on $B(r_{{\eta\psi^{(t)}}};y)$ and is thus lower bounded by integral of $\lambda_L$ on $B(r_{{\eta\psi^{(t)}}};y)$.

Let $n_{m}^{(t)}$ denotes the total number of samples fall into $B (r_{m};y)$ at time $t$. For $u<0$, we have the following upper bound:
\begin{equation}\begin{split}
	&P\{\Psi ^{(t)}<\lfloor\eta{\psi^{(t)}} \rfloor\}\\
	=&P\{ n_{\lfloor\eta{\psi^{(t)}} \rfloor}^{(t)}<\lfloor\eta{\psi^{(t)}} \rfloor\}\\
	\leq	&P\{ n_{\eta{\psi^{(t)}}}^{(t)}<\eta{\psi^{(t)}}  \}\\
	=&\int_{(p_{\eta\psi^{(t)}}^{(1)},p_{\eta\psi^{(t)}}^{(2)},\dots,p_{\eta\psi^{(t)}}^{(t)})}\\
	&P\{n_{\eta{\psi^{(t)}}}^{(t)}<\eta{\psi^{(t)}} | {(p_{\eta\psi^{(t)}}^{(1)},p_{\eta\psi^{(t)}}^{(2)},\dots,p_{\eta\psi^{(t)}}^{(t)})}\}
	dP\{{(p_{\eta\psi^{(t)}}^{(1)},p_{\eta\psi^{(t)}}^{(2)},\dots,p_{\eta\psi^{(t)}}^{(t)})}\}
\end{split}\label{inequality: conditional}
\end{equation}
Focusing on the conditional probability:
\begin{equation}
\begin{split}
	&P\{n_{\eta{\psi^{(t)}}}^{(t)}<\eta{\psi^{(t)}} | {(p_{\eta\psi^{(t)}}^{(1)},p_{\eta\psi^{(t)}}^{(2)},\dots,p_{\eta\psi^{(t)}}^{(t)})}\}\\
	\leq	&\inf_{u<0}\exp(-u\eta \psi^{(t)}+(e^u-1) \sum_{t'=1}^t  p_{\eta \psi^{(t)}}^{(t')}  )\\
	=&\exp(-\eta \psi^{(t)}\ln(\frac{\eta \psi^{(t)}}{ \sum_{t'=1}^t  p_{\eta \psi^{(t)}}^{(t')} })+\eta \psi^{(t)}- \sum_{t'=1}^t  p_{\eta \psi^{(t)}}^{(t')}  )\\
	\leq&\exp((-\eta ^{1+\xi}\ln\eta ^{1+\xi}+\eta ^{1+\xi}- 1 )\sum_{t'=1}^t  p_{\eta \psi^{(t)}}^{(t')} )\\
	\leq& \exp((-\eta ^{1+\xi}\ln\eta ^{1+\xi}+\eta ^{1+\xi}- 1 ) C(\eta \psi^{(t)})^{-\xi}t\lambda_L)\\
	=& \exp((-\eta ^{1+\xi}\ln\eta ^{1+\xi}+\eta ^{1+\xi}- 1 ) C(\eta (C\lambda_L)^\frac{1}{1+\xi}t^\frac{1}{1+\xi})^{-\xi}t\lambda_L)\\
	=&\exp(\eta^{-\xi}(-\eta ^{1+\xi}\ln\eta ^{1+\xi}+\eta ^{1+\xi}- 1 ) (C\lambda_L)^\frac{1}{1+\xi}t^\frac{1}{1+\xi}).
\end{split}\label{inequality: conditional probability}
\end{equation}
Based on (\ref{inequality: conditional}) and (\ref{inequality: conditional probability}), we have $\forall t\geq T^*$:
\begin{equation}
P\{\Psi ^{(t)}<\lfloor\eta{\psi^{(t)}} \rfloor\}\leq \exp(\eta^{-\xi}(-\eta ^{1+\xi}\ln\eta ^{1+\xi}+\eta ^{1+\xi}- 1 ) (C\lambda_L)^\frac{1}{1+\xi}t^\frac{1}{1+\xi})\label{inequation: insufficient observations}
\end{equation}
This leads to:
\begin{equation}
\begin{split}
	\lim_{t\to\infty}\sum_{t'=T^*}^tP\{\Psi ^{(t)}<\lfloor\eta{\psi^{(t)}} \rfloor\}\leq&\lim_{t\to\infty}\sum_{t'=T^*}^t\exp(\eta^{-\xi}(-\eta ^{1+\xi}\ln\eta ^{1+\xi}+\eta ^{1+\xi}- 1 ) (C\lambda_L)^\frac{1}{1+\xi}{t'}^\frac{1}{1+\xi})<\infty.
\end{split}\label{equation:1.1.3}
\end{equation}
The first inequality holds as the conditional probability share the same upper bound, thus the unconditional probability. The second inequality holds as $\eta^{-\xi}(-\eta ^{1+\xi}\ln\eta ^{1+\xi}+\eta ^{1+\xi}- 1 ) (C\lambda_L)^\frac{1}{1+\xi}$ is a negative constant and the series converges. Thus:
\begin{equation}
\sum_{t=1}^\infty P\{\Psi^{(t)} <\lfloor\eta\psi^{(t)} \rfloor\}<\infty.
\end{equation}Based on Borel-Cantelli lemma, this leads to:
\begin{equation}
P\{\limsup_{t\to\infty}\Psi^{(t)}<\lfloor\eta\psi^{(t)} \rfloor\}=0,
\end{equation}which finally leads to:
\begin{equation}
P\{\liminf_{t\to\infty}\Psi^{(t)}\geq\lfloor\eta\psi^{(t)} \rfloor\}=1.
\end{equation}
\endproof{$\blacksquare$}

\subsection{Lemma \ref{Lemma:error}}
\proof{Proof}

Since the discussion is for each $(x, y)$, we omit its notation when no confusion is caused. Consider the upper bound for the following probability for any arbitrary $\delta>0$:
\begin{equation}
\begin{split}
	&P\{|\bar{\epsilon}^{(t)}|>\delta, \Psi^{(t)} \geq\lfloor \eta {\psi^{(t)}} \rfloor\}\\
	\leq&\sup_{m\in[\lfloor \eta {\psi^{(t)}} \rfloor,t] , (\sigma_{[i_1 ]}, \dots, \sigma_{[i _{m }]})}P\{|\bar{\epsilon}^{(t)}|>\delta\Big| (\sigma_{[i_1 ]}, \dots, \sigma_{[i _{\Psi^{(t)} }]}) , {\Psi}^{(t)}=m \},
\end{split}\label{equation:probbound}
\end{equation}
where $i$ is relabeled as $i_j$ to index observations that fall into the neighborhood of our specified $(x,y)$. This is the probability such that when $t$ observations are accumulated, we not only have sufficient observations in the neighborhood (e.g.$\Psi^{(t)}\geq\lfloor\eta{\psi^{(t)}}\rfloor$) but also have error caused by noise greater than $\delta$. Notice, by definition:
\begin{equation}
\begin{split}
	\bar{\epsilon}^{(t)}  =\frac{\sum_{j=1}^{{\Psi^{(t)} }}\epsilon_{[i_j ]}}{{\Psi^{(t)}}},
\end{split}\label{equation:error upperbound}
\end{equation}
and conditioned on $(\sigma_{[i_1 ]}, \dots, \sigma_{[i _{\Psi^{(t)} }]})$, $\epsilon_{[i_j ]}$ follows $\mathcal{N}(0,\sigma_{[i_j ]})$ independently. Therefore, each conditional probability in $(\ref{equation:probbound})$ can be further bounded by:
\begin{equation}
\begin{split}
	&P\{|\bar{\epsilon}^{(t)}|>\delta\Big| (\sigma_{[i_1 ]}, \dots, \sigma_{[i _{\Psi^{(t)} }]}) , {\Psi}^{(t)}=m \}\\
	=&2\Phi(\frac{m(-\delta)}{\sqrt{\sum_{j=1}^{m}\sigma_{[i_j ]}^2}})\\
	\leq&2\Phi(\frac{m(-\delta)}{\sqrt{m(\sigma+L_2cm^{-\frac{\xi}{d}})^2}})\\
	\leq&2\exp(-\frac{m\delta^2}{2(\sigma+L_2cm^{-\frac{\xi}{d}})^2})\\
\end{split}
\label{equation:1.1}
\end{equation}
The first inequality holds due to the Lipschitz continuity of $\sigma_{[i_j ]}$ in the $Y$ neighborhood. The second inequality holds by Chernoff bound.	Therefore:
\begin{equation}
\begin{split}
	&\sup_{m\in[\lfloor \eta {\psi^{(t)}} \rfloor,t] , (\sigma_{[i_1 ]}, \dots, \sigma_{[i _{m }]})}P\{|\bar{\epsilon}^{(t)}|>\delta\Big| (\sigma_{[i_1 ]}, \dots, \sigma_{[i _{\Psi^{(t)} }]}) , {\Psi}^{(t)}=m \}\\
	\leq&\sup_{m\in[\lfloor \eta {\psi^{(t)}} \rfloor,t] , (\sigma_{[i_1 ]}, \dots, \sigma_{[i _{m }]})}2\exp(-\frac{m\delta^2}{2(\sigma+L_2cm^{-\frac{\xi}{d}})^2})\\
	=&\exp(-\frac{\lfloor \eta {\psi^{(t)}} \rfloor\delta^2}{2(\sigma+L_2c\lfloor \eta {\psi^{(t)}} \rfloor^{-\frac{\xi}{d}})^2}+\ln 2)\\
	\leq&\exp(-\frac{\lfloor \eta(C\lambda_L)^\frac{1}{1+\xi}t^\frac{1}{1+\xi}\rfloor\delta^2}{2(\sigma+L_2c\lfloor \eta(C\lambda_L)^\frac{1}{1+\xi}t^\frac{1}{1+\xi}\rfloor^{-\frac{\xi}{d}})^2}+\ln 2)\\
	\leq&\exp(-\frac{ [\eta (C\lambda_L)^\frac{1}{1+\xi}t^\frac{1}{1+\xi}-1]\delta^2}{2(\sigma+L_2c[\eta (C\lambda_L)^\frac{1}{1+\xi}t^\frac{1}{1+\xi}-1]^{-\frac{\xi}{d}})^2}+\ln 2)\\
	=&\exp(-t^\frac{1}{1+\xi}(\frac{ [\eta (C\lambda_L)^\frac{1}{1+\xi}-\frac{1}{t^\frac{1}{1+\xi}}]\delta^2}{2(\sigma+L_2c[\eta (C\lambda_L)^\frac{1}{1+\xi}t^\frac{1}{1+\xi}-1]^{-\frac{\xi}{d}})^2}+\frac{\ln 2}{t^{\frac{1}{1+\xi}}}))\\
	=:&\exp(-t^\frac{1}{1+\xi}H _{2}^{(t)}(\eta,\delta)),
\end{split}
\end{equation}
where $H _{2}^{(t)}(\eta,\delta):=(\frac{ [\eta (C\lambda_L)^\frac{1}{1+\xi}-\frac{1}{t^\frac{1}{1+\xi}}]\delta^2}{2(\sigma+L_2c[\eta (C\lambda_L)^\frac{1}{1+\xi}t^\frac{1}{1+\xi}-1]^{-\frac{\xi}{d}})^2}+\frac{\ln 2}{t^{\frac{1}{1+\xi}}})$. We have:
\begin{equation*}
H _{2}(\eta,\delta):=\lim_{t\to\infty}H _{2}^{(t)}(\eta,\delta)=\frac{ \eta (C\lambda_L)^\frac{1}{1+\xi}\delta^2}{2\sigma^2}>0,
\end{equation*}
where $\sigma$ is the abbreviation for $\sigma_x(y)$. This indicates that given $\eta\in(0,1)$, $\forall \delta>0, h_2\in(0, H_2(\eta,\delta))$, $\exists T^{**}$ such that $\forall t>T^{**}$, $H_{2}^{(t)}(\eta,\delta)>h_2$.  Based on (\ref{equation:probbound}), we have $\forall t>T^{**}$:
\begin{equation}
P\{|\bar{\epsilon}^{(t)}|>\delta, \Psi^{(t)}  \geq\lfloor \eta {\psi^{(t)}} \rfloor\}\leq\exp[-t^{\frac{1}{1+\xi}}H_{2}^{(t)} ]\leq \exp[-t^{\frac{1}{1+\xi}}h_2 ]\label{equation:secondterm}
\end{equation}
Therefore:
\begin{equation}
\begin{split}
	&\lim_{t\to\infty}\sum_{t'=T^{**}}^tP\{|\bar{\epsilon}^{(t')}|>\delta, \Psi ^{(t')}  \geq\lfloor \eta {\psi^{(t')}} \rfloor\}\\
	\leq&\lim_{t\to\infty}\sum_{t'=T^{**}}^t\exp[-{t'}^{\frac{1}{1+\xi}}h_2 ]<\infty
\end{split}
\label{equation:lowerbound2}
\end{equation}
This means the series of sum of probability is bounded $\forall \delta>0$:
$$\lim_{t\to\infty}\sum_{t'=1}^tP\{|\bar{\epsilon}^{(t')}|>\delta, \Psi^{(t')}  \geq\lfloor \eta {\psi^{(t')}} \rfloor\}<\infty$$
By Borel-Cantelli Lemma, we can conclude that:
$$P\{\lim_{t\to\infty}|\bar{\epsilon}^{(t)}|=0\cup\Psi ^{(t)}<\lfloor \eta {\psi^{(t)}} \rfloor\}=1
$$ 
Since by Lemma \ref{Lemma:mt}, $P\{\lim_{t\to\infty}\Psi ^{(t)}\geq\lfloor \eta {\psi^{(t)}} \rfloor\}=1$, we have:
\begin{equation}\begin{split}
	1=&P\{(\lim_{t\to\infty}|\bar{\epsilon}^{(t)}   |=0\cup\Psi ^{(t)}<\lfloor \eta {\psi^{(t)}} \rfloor)\cap(\lim_{t\to\infty}\Psi ^{(t)}\geq\lfloor \eta {\psi^{(t)}} \rfloor)\}\\
	=&P\{\lim_{t\to\infty}|\bar{\epsilon}^{(t)}   |=0\cap\Psi ^{(t)}\geq\lfloor \eta {\psi^{(t)}} \rfloor\}\\
	\leq&P\{\lim_{t\to\infty}|\bar{\epsilon}^{(t)}   |=0\}=1
\end{split}
\end{equation}
Thus complete the prove.
\endproof{$\blacksquare$}

\subsection{Theorem \ref{theorem: 1st 2nd}}
\proof{Proof}

Since the discussion is for each $(x, y)$, we omit its notation when no confusion is caused (e.g. $f$, $\sigma$ indicating $f_x(y)$ and $\sigma_x(y)$). We relabel $i$ using $i_j$ to label the points fall into the neighborhood of $(x,y)$. 
\subsubsection{Consistency of the first order information} 
We have:
\begin{equation}
\begin{split}
	{\hat{f}}^{(t)}-f=&\frac{1}{{\Psi^{(t)}}}\sum_{j=1}^{\Psi^{(t)}} \tilde{G}_{[i_j]}-f_{[i_j]}+f_{[i_j]}-f\\
	=&\bar{\epsilon}^{(t)} +\frac{1}{{\Psi^{(t)}}}\sum_{j=1}^{\Psi^{(t)}} [f_{[i_j]}-f]
\end{split}\label{equation:1}
\end{equation}
The aforementioned formulation separates the estimation error into two terms with the first term $\bar{\epsilon}^{(t)} $ representing error caused by noise and the second term $\frac{1}{{\Psi^{(t)}}}\sum_{j=1}^{\Psi^{(t)}} [f_{[i_j]}-f]$ representing error caused by bias.

For $\bar{\epsilon}^{(t)}   $, by Lemma \ref{Lemma:error}, we proved that $\bar{\epsilon}^{(t)}\to 0$ with probability 1. For the bias term, since we assume that the function $f$ is Lipschitz continuous with constant $L_1$, it is bounded by the size of the ball:
\begin{equation*}
|\frac{1}{{\Psi^{(t)}}}\sum_{j=1}^{\Psi^{(t)}} [f_{[i_j]}-f]|\leq{r_{\Psi ^{(t)}}}L_1
\end{equation*}
As $\lim_{t\to\infty} {\psi^{(t)}} =\infty$ and by Lemma \ref{Lemma:mt}, $\Psi^{(t)}  \geq\lfloor \eta {\psi^{(t)}} \rfloor$ with probability 1. These imply that when $t\to\infty$:
\begin{equation}
|\frac{1}{{\Psi^{(t)}}}\sum_{j=1}^{\Psi^{(t)}} [f_{[i_j]}-f]|\leq{r_{\Psi ^{(t)}}}L_1 \overset{a.s.}{\leq}{r_{ \eta {\psi^{(t)}}  }}L_1\to0, \label{equation:3}
\end{equation}
meaning $$	|\frac{1}{{\Psi^{(t)}}}\sum_{j=1}^{\Psi^{(t)}} [f_{[i_j]}-f]|\to 0\quad a.s.$$ By Continuous Mapping Theorem, when $t\to\infty$, $${\hat{f}}^{(t)}-f=\bar{\epsilon}^{(t)}+\frac{1}{{\Psi^{(t)}}}\sum_{j=1}^{\Psi^{(t)}} [f_{[i_j]}-f]{\to}0\quad a.s.$$
\subsubsection{Consistency of the second order information} 
From definition of sample variance within the neighborhood:
\begin{equation}
{[\hat{\sigma}^{(t)}]}^2
=[\frac{1}{{\Psi^{(t)}}-1}\sum_{j=1}^{\Psi^{(t)}} (\tilde{G}_{[i_j]}-f )^2]-[\frac{\Psi^{(t)}}{\Psi^{(t)} -1}({\hat{f}^{(t)}  }-f )^2]\label{definition:sample variance}
\end{equation}
For the second term, based on Theorem \ref{theorem: 1st 2nd}, ${\hat{f}^{(t)}  }-f \to0$ $a.s.$ In addition $\frac{\Psi^{(t)}}{\Psi^{(t)} -1}\to 0$ $a.s.$  By Continuous Mapping Theorem, 
\begin{equation}
\frac{\Psi^{(t)}}{\Psi^{(t)} -1}({\hat{f}^{(t)}  }-f )^2\to 0\quad a.s.\label{equation:3.1}
\end{equation}
For the first term of equation (\ref{definition:sample variance}): 
\begin{equation}
\begin{split}
	&\frac{1}{{\Psi^{(t)}}-1}\sum_{j=1}^{\Psi^{(t)}} (\tilde{G}_{[i_j]}-f )^2\\
	=&\frac{\Psi^{(t)}}{\Psi^{(t)}-1}\sum_{j=1}^{\Psi^{(t)}} \frac{1}{\Psi^{(t)}}\epsilon_{[i_j ]}^2+\frac{1}{\Psi^{(t)}}(f_{[i_j]}-f )^2+\frac{2}{\Psi^{(t)}}\epsilon_{[i_j]}(f_{[i_j]}-f )
\end{split}\label{equation:second order split}
\end{equation}
Now we analyze each term in (\ref{equation:second order split}). Since $\frac{\Psi^{(t)}}{\Psi^{(t)}-1}\to 1$ $a.s.$, we can ignore it for the moment. For $\sum_{j=1}^{\Psi^{(t)}} \frac{1}{\Psi^{(t)}}\epsilon_{[i_j ]}^2$, due to Lipschitz continuity, it is bounded by:
\begin{equation}
\sum_{j=1}^{\Psi^{(t)}} \frac{1}{\Psi^{(t)}}\epsilon_{[i_j ]}^2
\leq({\sigma }+L_2r_{{\Psi ^{(t)}  }})^2\frac{\sum_{j=1}^{\Psi^{(t)}}  (z_{[i_j]})^2}{{\Psi^{(t)} }},\label{equation:2.0}
\end{equation}
where $z_{[i_j]}$ is the standardized normal variable for each $\epsilon_{[i_j]}$ (i.e. $z_{[i_j]}:=\epsilon_{[i_j]}/\sigma_{[i_j]}$).
We now show that $\frac{\sum_{j=1}^{\Psi^{(t)}}  (z_{[i_j]})^2}{{\Psi^{(t)} }}\to 1$ with probability 1 using Borel-Cantelli Lemma. For an arbitrary $\delta>0$:
\begin{equation}
\begin{split}
	&P\{|\frac{\sum_{j=1}^{\Psi^{(t)}}  (z_{[i_j]})^2}{{\Psi^{(t)} }}-1|>\delta\}\\
	\leq&P\{|\frac{\sum_{j=1}^{\Psi^{(t)}}  (z_{[i_j]})^2}{{\Psi^{(t)} }}-1|>\delta, \Psi^{(t)}  \geq\lfloor \eta {\psi^{(t)}} \rfloor\}\\
	&+P\{\Psi^{(t)} <\lfloor \eta {\psi^{(t)}} \rfloor\}
\end{split}\label{equation:2.1}
\end{equation}
For the first part of equation (\ref{equation:2.1}), we have that:
\begin{align*}
&P\{|\frac{\sum_{j=1}^{\Psi^{(t)}}  (z_{[i_j]})^2}{{\Psi^{(t)} }}-1|>\delta, \Psi^{(t)}  \geq\lfloor \eta {\psi^{(t)}} \rfloor\}\\
\leq&\max_{m\geq \lfloor{ \eta {\psi^{(t)}} }\rfloor} P\{|{\sum_{j=1}^{\Psi^{(t)}}  (z_{[i_j]})^2}-{\Psi^{(t)} }|>\delta{\Psi^{(t)} }\Big| \Psi^{(t)}=m\}\\
=&\max_{m\geq \lfloor{ \eta {\psi^{(t)}} }\rfloor} P\{{\sum_{j=1}^{\Psi^{(t)}}  (z_{[i_j]})^2}>(1+\delta){\Psi^{(t)} }\Big| \Psi^{(t)}=m\}+P\{{\sum_{j=1}^{\Psi^{(t)}}  (z_{[i_j]})^2}<(1-\delta){\Psi^{(t)} }\Big| \Psi^{(t)}=m\}
\end{align*}
Conditioned on $\Psi^{(t)}=m$,  ${\sum_{j=1}^{\Psi^{(t)}}  (z_{[i_j]})^2}$ follows $\chi^2$ distribution with $m$ degree of freedom. By using Chernoff Bound, we have that conditioned on each $m$:
\begin{equation*}
\begin{split}
	&P\{{\sum_{j=1}^{\Psi^{(t)}}  (z_{[i_j]})^2}>(1+\delta){\Psi^{(t)} }\Big| \Psi^{(t)}=m\}+P\{{\sum_{j=1}^{\Psi^{(t)}}  (z_{[i_j]})^2}<(1-\delta){\Psi^{(t)} }\Big| \Psi^{(t)}=m\}\\
	\leq&((1+\delta)e^{\delta})^\frac{m}{2}+((1-\delta)e^{\delta})^\frac{m}{2}
\end{split}
\end{equation*}
For a small $\delta>0$ such that $(1+\delta)e^{\delta}<1$, the value takes the maximum when $m$ is the smallest (i.e. $m=\lfloor{ \eta {\psi^{(t)}} }\rfloor$). For such $\delta$, we have that:
\begin{align*}
&P\{|\frac{\sum_{j=1}^{\Psi^{(t)}}  (z_{[i_j]})^2}{{\Psi^{(t)} }}-1|>\delta, \Psi^{(t)}  \geq\lfloor \eta {\psi^{(t)}} \rfloor\}\\
\leq&((1+\delta)e^{\delta})^\frac{\lfloor{ \eta {\psi^{(t)}} }\rfloor}{2}+((1-\delta)e^{\delta})^\frac{\lfloor{ \eta {\psi^{(t)}} }\rfloor}{2}\\
=&((1+\delta)e^{\delta})^\frac{\lfloor{ \eta (C\lambda_L)^\frac{1}{1+\xi}t^\frac{1}{1+\xi} }\rfloor}{2}+((1-\delta)e^{\delta})^\frac{\lfloor{ \eta (C\lambda_L)^\frac{1}{1+\xi}t^\frac{1}{1+\xi} }\rfloor}{2}
\end{align*} 
Since $\xi\in(0, 1)$, we conclude that:
\begin{equation}
\begin{split}
	&\lim_{t'\to\infty}\sum_{t=1}^{t'}P\{|\frac{\sum_{j=1}^{\Psi^{(t)}}  (z_{[i_j]})^2}{{\Psi^{(t)} }}-1|>\delta, \Psi^{(t)}  \geq\lfloor \eta {\psi^{(t)}} \rfloor\}\\
	\leq &\lim_{t'\to\infty}\sum_{t=1}^{t'}((1+\delta)e^{\delta})^\frac{\lfloor{ \eta (C\lambda_L)^\frac{1}{1+\xi}t^\frac{1}{1+\xi} }\rfloor}{2}+((1-\delta)e^{\delta})^\frac{\lfloor{ \eta (C\lambda_L)^\frac{1}{1+\xi}t^\frac{1}{1+\xi} }\rfloor}{2}<\infty
\end{split}
\label{equation:3.2}
\end{equation} 
For the second part of equation (\ref{equation:2.1}), based on Lemma \ref{Lemma:mt}, 
\begin{equation}
\lim_{t'\to\infty}\sum_{t=1}^{t} P\{\Psi <\lfloor \eta {\psi^{(t)}} \rfloor\}<\infty\label{equation:3.3}
\end{equation}
By combining equation (\ref{equation:3.2}) and (\ref{equation:3.3}), we get
\begin{equation}
\lim_{t'\to\infty}\sum_{t=1}^{t'}P\{|\frac{\sum_{j=1}^{\Psi^{(t)}}  (z_{[i_j]})^2}{{\Psi^{(t)} }}-1|>\delta\}<\infty
\end{equation}
By Borel-Cantelli lemma:
\begin{equation}
\frac{\sum_{j=1}^{\Psi^{(t)}}  (z_{[i_j]})^2}{{\Psi^{(t)} }}\to 0\quad a.s.\label{equation:3.4}
\end{equation}
Combined with the fact that $r_{{\Psi ^{(t)}  }}\to 0$ with probability 1, based on Continuous Mapping Theorem, the equation (\ref{equation:2.0}) shows that:
\begin{equation*}
\begin{split}
	[{\sigma }]^2\overset{a.s.}{\leftarrow}&({\sigma }-L_2r_{{\Psi   }})^2\frac{\sum_{j=1}^{\Psi^{(t)}}  (z_{[i_j]})^2}{{\Psi^{(t)} }}\\
	\leq&\sum_{j=1}^{\Psi^{(t)}} \frac{1}{\Psi^{(t)}}\epsilon_{[i_j ]}^2\\
	\leq&({\sigma }+L_2r_{{\Psi   }})^2\frac{\sum_{j=1}^{\Psi^{(t)}}  (z_{[i_j]})^2}{{\Psi^{(t)} }}\\
	\overset{a.s.}{\to}& [{\sigma }]^2
\end{split}
\end{equation*}
This suggests that
\begin{equation}
\sum_{j=1}^{\Psi^{(t)}} \frac{1}{\Psi^{(t)}}\epsilon_{[i_j ]}^2{\to}[{\sigma }]^2\quad a.s.\label{equation:2.9.1}
\end{equation} 
For the second part in equation (\ref{equation:second order split}), since from Lemma 2, we can conclude that ${\Psi^{(t)}   }\to\infty$ with probability 1. By Continuous Mapping Theorem, $r_{{\Psi ^{(t)}  }}\to 0$ with probability 1. Therefore, we have:
\begin{align*}
0\leq\frac{1}{\Psi^{(t)}}\sum_{j=1}^{\Psi^{(t)}} (f_{[i_j]}-f )^2\leq(Lr_{{\Psi   }})^2\overset{a.s.}{\to}0
\end{align*}
This suggests that:
\begin{equation}
\frac{1}{\Psi^{(t)}}\sum_{j=1}^{\Psi^{(t)}} (f_{[i_j]}-f )^2\to 0\quad a.s. \label{equation:3.5}
\end{equation}
For the third part in equation (\ref{equation:second order split}), it is bounded by:
\begin{equation*}
-2L_1r_{{\Psi ^{(t)}  }}        \bar{\epsilon} ^{(t)}       \leq\frac{2}{\Psi^{(t)}}\sum_{j=1}^{\Psi^{(t)}} \epsilon_{[i_j]}(f_{[i_j]}-f )\leq2L_1r_{{\Psi ^{(t)}  }}        \bar{\epsilon}^{(t)}    
\end{equation*}
From Lemma \ref{Lemma:error}, $\bar{\epsilon}^{(t)}    \to 0$ with probability 1. Therefore, we have: 
\begin{equation*}
\begin{split}
	0\overset{a.s.}{\leftarrow}&	-2L_1r_{{\Psi ^{(t)}  }}        \bar{\epsilon} ^{(t)}   \\
	\leq&\frac{2}{\Psi^{(t)}}\sum_{j=1}^{\Psi^{(t)}} \epsilon_{[i_j]}(f_{[i_j]}-f )\\
	\leq&2L_1r_{{\Psi ^{(t)}  }}        \bar{\epsilon}^{(t)}      \overset{a.s.}{\to}0
\end{split}
\end{equation*}
This implies that:
\begin{equation}
\frac{2}{\Psi^{(t)}}\sum_{j=1}^{\Psi^{(t)}} \epsilon_{[i_j]}(f_{[i_j]}-f )\to0\quad a.s.\label{equation:3.6}
\end{equation}
Based on formulation (\ref{equation:3.1}) (\ref{equation:2.9.1}) (\ref{equation:3.5}) and (\ref{equation:3.6}), by Continuous Mapping Theorem, the statements is proved.

\endproof{$\blacksquare$}

\subsection{Proposition \ref{proposition: consistency}}
\proof{Proof}Based on the definition, the total decision loss is formulated as:
\begin{equation*}
\varpi^{(t)} =\int_{y\in Y^\circ}\upsilon^{(t)}(y)dP_Y=\int_{y\in Y^\circ}E_{\mathcal{D}^{(t)} }d({f^\dagger}^{(t)}(y),{f^*(y)}) dP_Y.
\end{equation*}
We first show that,$\forall y\in Y^\circ$, $d({f^\dagger}^{(t)}(y),{f^*(y)})\overset{a.s.}{\to}0$.
By Theorem \ref{theorem: 1st 2nd}, for any $x$ and $y$, ${\hat{f}_x^{(t)}(y)}\overset{a.s.}{\to}f_{x}(y)$ when $t\to\infty$. Since $\forall x\notin x^*(y)$, $f_x(y)>f_{x^*}(y)$, we have that when $t\to\infty$, 
$$\min_{x\notin x^*(y)}{\hat{f}_x^{(t)}(y)}> \hat{f}_{x^*}^{(t)}(y)\quad a.s.$$ 
This suggests that when $t\to\infty$, $\forall y\in Y^\circ$, $\argmin_{x\in\chi}   \hat{f}_x^{(t)}(y) =:{\hat{x}^{*(t)}}(y)   \in x^*(y), a.s.$ 	 Therefore:  
$$d({f^\dagger}^{(t)}(y),{f^*(y)})\to0\quad a.s.$$
As the almost surely convergence implies convergence in expectation, this leads to the fact that $$\lim_{t\to\infty} \upsilon^{(t)} (y)=\lim_{t\to\infty} E_{\mathcal{D}^{(t)}}  d({f^\dagger}^{(t)}(y),{f^*(y)}) =0,$$ $\forall y\in Y^\circ$. Notice that $f_x(y)$ is uniformly bounded on $Y^\circ$ by Lipschitz continuity assumed in Assumption \ref{assumption: lipschitz} and the fact that $Y$ is bounded. This suggests that $d(f^\dagger(y),{f^*(y)})$ is uniformly bounded on $Y^\circ$.  Based on Dominated Convergence Theorem, the interchange of limit and the integral is allowed: 
\begin{equation*}
\begin{split}
	\lim_{t\to\infty}\varpi ^{(t)}=&\lim_{t\to\infty}\int_{y\in Y^\circ} \upsilon^{(t)} (y)dP_Y(y)\\
	=&\int_{y\in Y^\circ}\lim_{t\to\infty} \upsilon^{(t)} (y)dP_Y(y)\\
	=0
\end{split}
\end{equation*}
Thus complete the proof.
\endproof{$\blacksquare$}

\subsection{Theorem \ref{theorem: CDL convergence rate}}
\proof{Proof}
$\forall y\in Y\backslash MO$, $\upsilon^{(t)} $ is upper bounded by:
\begin{equation}
\begin{split}
	\upsilon^{(t)} (y)=&E_{\mathcal{D}^{(t)} }d({f^\dagger}^{(t)}(y),{f^*(y)})\\
	=&\sum_{x\neq x^*(y)}d({f^\dagger}^{(t)}(y),{f^*(y)})P\{\hat{x}^{*(t)}(y)=x\}\\
	\leq&(|\chi|-1)\max_{x\neq x^*(y)}d({f^\dagger}^{(t)}(y),{f^*(y)})P\{{\hat{f}_x^{(t)}(y)}<\hat{f}_{x^*}^{(t)}(y)\}
\end{split}\label{equation:CDLupper}
\end{equation}
Now we focus on $P\{{\hat{f}_x^{(t)}(y)}<\hat{f}_{x^*}^{(t)}(y)\}$. For each $x\neq x^*(y)$, this probability can be separated into the one caused by observation noise and the one caused by bias. We first show the upper bound of the probability:
\begin{align*}
&P\{\hat{f}_{x}^{(t)}(y)<\hat{f}_{x^*}^{(t)}(y)\}\\
=&P\{\bar{\epsilon}^{(t)} _{x, y}+ \frac{1}{{\Psi^{(t)}_{x,y}}}\sum_{i:x_{[i]}=x, y_{[i]}\in {b _{x, y}}}[f_{[i]}-f_{x}(y)]+f_{x}(y)   \\
&<{\bar{\epsilon}^{(t)} _{x^*, y}}+ \frac{1}{{\Psi_{x^*, y}^{(t)}}}\sum_{i:x_{[i]}= x^*(y), y_{[i]}=b _{x^*, y}}[f_{[i]}-{f_{x^*}(y)}]+{f_{x^*}(y)}  \}\\
\leq&P\{{\bar{\epsilon}^{(t)} _{x^*, y}}-\bar{\epsilon}^{(t)} _{x, y}>{\delta_{x, y}^{(t)}}\}, 
\end{align*}
where $\forall x\neq x^*(y)$, ${\delta_{x, y}^{(t)}}={\delta_{x, y}^{(t)}}(\Psi_{x^*, y}^{(t)},\Psi_{x, y}^{(t)}):=-({r_{\Psi_{x^*, y}^{(t)}}}+{r_{\Psi_{x, y}^{(t)}}})L_1  + f_{x}(y)-{f_{x^*}(y)}$. The inequality holds due to Lipschitz continuity of $f_x(y)$.
For all cases where $\Psi_{x, y}^{(t)}<\lfloor\eta\psi^{(t)}\rfloor$ for some $\eta\in(0,1)$, they are classified as having insufficient observations in the neighborhood of $(x,y)$. Therefore:
\begin{equation}
\begin{split}
	&P\{{\bar{\epsilon}^{(t)} _{x^*, y}}-\bar{\epsilon}^{(t)} _{x, y}> {\delta_{x, y}^{(t)}} \}    \\
	=&P\{{\bar{\epsilon}^{(t)} _{x^*, y}}-\bar{\epsilon}^{(t)} _{x, y}> {\delta_{x, y}^{(t)}} , \Psi_{x, y}^{(t)}< \lfloor\eta\psi ^{(t)}\rfloor \cup\Psi_{x^*, y}^{(t)}< \lfloor\eta\psi ^{(t)}\rfloor\}    \\
	&+P\{{\bar{\epsilon}^{(t)} _{x^*, y}}-\bar{\epsilon}^{(t)} _{x, y}> {\delta_{x, y}^{(t)}} , \Psi_{x, y}^{(t)}\geq \lfloor\eta\psi ^{(t)}\rfloor, \Psi_{x^*, y}^{(t)}\geq \lfloor\eta\psi ^{(t)}\rfloor\}    \\
	\leq&P\{\Psi_{x, y}^{(t)}< \lfloor\eta\psi ^{(t)}\rfloor \}+P\{\Psi_{x^*, y}^{(t)}< \lfloor\eta\psi ^{(t)}\rfloor\}\\
	&+P\{{\bar{\epsilon}^{(t)} _{x^*, y}}-\bar{\epsilon}^{(t)} _{x, y}> {\delta_{x, y}^{(t)}}, \Psi_{x, y}^{(t)}\geq \lfloor\eta\psi ^{(t)}\rfloor , \Psi_{x^*, y}^{(t)}\geq \lfloor\eta\psi ^{(t)}\rfloor\}
\end{split}\label{equation:separation}
\end{equation}
For the first two terms, based on formulation (\ref{inequation: insufficient observations}), we have that exists $T^*$, such that $\forall t\geq T^*$:
\begin{equation}
\begin{split}
	&P\{\Psi_{x, y}^{(t)}< \lfloor\eta\psi ^{(t)}\rfloor \}+P\{\Psi_{x^*, y}^{(t)}< \lfloor\eta\psi ^{(t)}\rfloor\}\\
	\leq&2\exp(\eta^{-\xi}(-\eta ^{1+\xi}\ln\eta ^{1+\xi}+\eta ^{1+\xi}- 1 ) (C\lambda_L)^\frac{1}{1+\xi}t^\frac{1}{1+\xi}).
\end{split}\label{equation: errortwo}
\end{equation}
For the third term:
\begin{equation}
\begin{split}
	&P\{{\bar{\epsilon}^{(t)} _{x^*, y}}-\bar{\epsilon}^{(t)} _{x, y}> \delta_{x, y}^{(t)} , \Psi_{x, y}^{(t)}\geq \lfloor\eta\psi ^{(t)}\rfloor , \Psi_{x^*, y}^{(t)}\geq \lfloor\eta\psi ^{(t)}\rfloor\}  \\
	\leq&\sup_{\Psi_{x^*, y}^{(t)}\geq \lfloor\eta\psi ^{(t)}\rfloor, \Psi_{x, y}^{(t)}\geq \lfloor\eta\psi ^{(t)}\rfloor,(\sigma_{[i_1]},\sigma_{[i_2]},\dots,\sigma_{[i_{\Psi_{x^*, y}^{(t)}}]}),(\sigma_{[k_1]},\sigma_{[k_2]},\dots,\sigma_{[k_{\Psi^{(t)}_{x,y}}]})}\\
	&P\{{\bar{\epsilon}^{(t)} _{x^*, y}}-\bar{\epsilon}^{(t)} _{x, y}> \delta_{x, y}^{(t)} \Big|(\sigma_{[i_1]},\sigma_{[i_2]},\dots,\sigma_{[i_{\Psi_{x^*, y}^{(t)}}]}),(\sigma_{[k_1]},\sigma_{[k_2]},\dots,\sigma_{[k_{\Psi^{(t)}_{x,y}}]}),\\
	&{\Psi_{x^*, y}^{(t)}}, {\Psi_{x, y}^{(t)} }\}\\
	=&\sup_{\Psi_{x^*, y}^{(t)}\geq \lfloor\eta\psi ^{(t)}\rfloor, \Psi_{x, y}^{(t)}\geq \lfloor\eta\psi ^{(t)}\rfloor,(\sigma_{[i_1]},\sigma_{[i_2]},\dots,\sigma_{[i_{\Psi_{x^*, y}^{(t)}}]}),(\sigma_{[k_1]},\sigma_{[k_2]},\dots,\sigma_{[k_{\Psi^{(t)}_{x,y}}]})}\\
	&\Phi(-\frac{\delta_{x, y}^{(t)}(\Psi_{x^*, y}^{(t)},\Psi_{x, y}^{(t)})}{\sqrt{\frac{\sum_{j=1}^{{\Psi_{x^*,y}}^{(t)}}\sigma_{x_{[i_j]}}^2(y_{[i_j]})}{({\Psi_{x^*,y}}^{(t)})^2}+\frac{\sum_{j=1}^{{\Psi_{x,y}}^{(t)}}\sigma_{x_{[i_k]}}^2(y_{[i_k]})}{({\Psi_{x,y}}^{(t)})^2}}})\\
	\leq&\Phi(-\frac{{\delta_{x, y}^{(t)}}(\lfloor\eta\psi ^{(t)}\rfloor,\lfloor\eta\psi ^{(t)}\rfloor)}
	{\sqrt{\frac{(\sigma_{x^*}(y)+L_2c\lfloor\eta\psi ^{(t)}\rfloor^{-\frac{\xi}{d}})^2}
			{\lfloor\eta\psi ^{(t)}\rfloor}
			+\frac{(\sigma_{x}(y)+L_2c\lfloor\eta\psi ^{(t)}\rfloor^{-\frac{\xi}{d}})^2}
			{\lfloor\eta\psi ^{(t)}\rfloor}}})\\
	\leq&\exp(-\frac{{\delta_{x, y}^{(t)}}^2(\lfloor\eta\psi ^{(t)}\rfloor,\lfloor\eta\psi ^{(t)}\rfloor)}
	{2(\frac{(\sigma_{x^*}(y)+L_2c\lfloor\eta\psi ^{(t)}\rfloor^{-\frac{\xi}{d}})^2}
		{\lfloor\eta\psi ^{(t)}\rfloor}
		+\frac{(\sigma_{x}(y)+L_2c\lfloor\eta\psi ^{(t)}\rfloor^{-\frac{\xi}{d}})^2}
		{\lfloor\eta\psi ^{(t)}\rfloor})})\\
	<&\exp(-\frac{{\delta_{x, y}^{(t)}}^2(\eta(C\lambda_L)^\frac{1}{1+\xi}t^\frac{1}{1+\xi}-1,\eta((C\lambda_L)^\frac{1}{1+\xi}t^\frac{1}{1+\xi}-1))}
	{2(\frac{(\sigma_{x^*}(y)+L_2c(\eta((C\lambda_L)^\frac{1}{1+\xi}t^\frac{1}{1+\xi}-1))^{-\frac{\xi}{d}})^2}
		{(\eta((C\lambda_L)^\frac{1}{1+\xi}t^\frac{1}{1+\xi}-1))}
		+\frac{(\sigma_{x}(y)+L_2c(\eta((C\lambda_L)^\frac{1}{1+\xi}t^\frac{1}{1+\xi}-1))^{-\frac{\xi}{d}})^2}
		{(\eta((C\lambda_L)^\frac{1}{1+\xi}t^\frac{1}{1+\xi}-1))})})\\
	\leq&\exp(-t^{\frac{1}{1+\xi}}\frac{{\delta_{x, y}^{(t)}}^2((\eta((C\lambda_L)^\frac{1}{1+\xi}t^\frac{1}{1+\xi}-1)),(\eta((C\lambda_L)^\frac{1}{1+\xi}t^\frac{1}{1+\xi}-1)))}
	{2(\frac{(\sigma_{x^*}(y)+L_2c(\eta((C\lambda_L)^\frac{1}{1+\xi}t^\frac{1}{1+\xi}-1))^{-\frac{\xi}{d}})^2}
		{\eta(C\lambda_L)^\frac{1}{1+\xi}-\frac{\eta+1}{t^\frac{1}{1+\xi}}}
		+\frac{(\sigma_{x}(y)+L_2c(\eta((C\lambda_L)^\frac{1}{1+\xi}t^\frac{1}{1+\xi}-1))^{-\frac{\xi}{d}})^2}
		{\eta(C\lambda_L)^\frac{1}{1+\xi}-\frac{\eta+1}{t^\frac{1}{1+\xi}}})})\\
	=:&\exp(-t^{\frac{1}{1+\xi}} H_3^{(t)}(x,y;\eta)),
\end{split}\label{equation:mixed probability upper bound}
\end{equation}
where $
H_3^{(t)}(x,y;\eta)
:=\frac{{\delta_{x, y}^{(t)}}^2((\eta((C\lambda_L)^\frac{1}{1+\xi}t^\frac{1}{1+\xi}-1)),(\eta((C\lambda_L)^\frac{1}{1+\xi}t^\frac{1}{1+\xi}-1)))}
{2(\frac{(\sigma_{x^*}(y)+L_2c(\eta((C\lambda_L)^\frac{1}{1+\xi}t^\frac{1}{1+\xi}-1))^{-\frac{\xi}{d}})^2}
{\eta(C\lambda_L)^\frac{1}{1+\xi}-\frac{\eta+1}{t^\frac{1}{1+\xi}}}
+\frac{(\sigma_{x}(y)+L_2c(\eta((C\lambda_L)^\frac{1}{1+\xi}t^\frac{1}{1+\xi}-1))^{-\frac{\xi}{d}})^2}
{\eta(C\lambda_L)^\frac{1}{1+\xi}-\frac{\eta+1}{t^\frac{1}{1+\xi}}})}.
$
The first inequality holds due to the maximum of conditional probability is no less than the whole probability. The first equality holds as conditioned on $\Psi_{x, y}^{(t)}$ and $\Psi_{x^*, y}^{(t)}$, $\delta_{x, y}^{(t)}$ is a constant. Plus the fact that conditioned on $\sigma_{[i]}$ for each noise term within neighborhood of $(x,y)$ and $(x^*,y)$, ${\bar{\epsilon}^{(t)} _{x^*, y}}-\bar{\epsilon}^{(t)} _{x, y}$ follows normal distribution. The second inequality holds due to Lipschitz continuity of the $\sigma_x(y)$. The third inequality holds due to Chernoff bound. The fourth inequality holds as $\lfloor\eta\psi^{(t)}\rfloor\geq\eta\psi^{(t)}-1$ and the function is again monotone decreasing w.r.t. to the entry. Therefore, based on (\ref{equation:separation}), (\ref{equation:separation}) and (\ref{equation:mixed probability upper bound}), we have:
\begin{equation}
P\{\hat{f}_{x}^{(t)}(y)<\hat{f}_{x^*}^{(t)}(y)\}\leq 2\exp(\eta^{-\xi}(-\eta ^{1+\xi}\ln\eta ^{1+\xi}+\eta ^{1+\xi}- 1 ) (C\lambda_L)^\frac{1}{1+\xi}t^\frac{1}{1+\xi})+ \exp(-t^{\frac{1}{1+\xi}} H_3^{(t)}(x,y;\eta))\label{inequality: mutual probability}
\end{equation}
$H_3^{(t)}(x,y;\eta)$ converges to a constant:
\begin{equation}
\lim_{t\to\infty}H_3^{(t)}(x,y;\eta)=\eta C^\frac{1}{1+\xi} \lambda_L^\frac{1}{1+\xi}\frac{(f_x(y)-f_{x^*}(y))^2}{2({\sigma_{x^*}^2(y)}+{\sigma_{x}^2(y)})}=:H_3(x,y;\eta)
\end{equation}
Therefore, we have:
\begin{equation}
\begin{split}
	&\liminf_{t\to\infty}-\frac{1}{t^{\frac{1}{1+\xi}}}\ln \upsilon^{(t)} (y)\\
	\geq&\liminf_{t\to\infty}-\frac{1}{t^{\frac{1}{1+\xi}}}\ln\Big[(|\chi|-1)\max_{x\neq x^*(y)}P\{{\hat{f}_x^{(t)}(y)}<\hat{f}_{x^*}^{(t)}(y)\}\max_{x\neq x^*(y)}d({f^\dagger}^{(t)}(y),{f^*(y)})\Big]\\
	=&\liminf_{t\to\infty}-\frac{1}{t^{\frac{1}{1+\xi}}}\ln \max_{x\neq x^*(y)} P\{\hat{f}_x^{(t)}(y)<\hat{f}_{x^*}^{(t)}(y)\}\\
	=&\liminf_{t\to\infty}\min_{x\neq x^*(y)}-\frac{1}{t^{\frac{1}{1+\xi}}}\ln  P\{\hat{f}_x^{(t)}(y)<\hat{f}_{x^*}^{(t)}(y)\}\\
	\geq &\liminf_{t\to\infty}\min_{x\neq x^*(y)}-\frac{1}{t^{\frac{1}{1+\xi}}}\ln[\exp(\eta^{-\xi}(-\eta ^{1+\xi}\ln\eta ^{1+\xi}+\eta ^{1+\xi}- 1 ) (C\lambda_L)^\frac{1}{1+\xi}t^\frac{1}{1+\xi}+\ln 2)\\
	&+\exp(-t^{\frac{1}{1+\xi}} H_3^{(t)}(x,y;\eta))]\\
	\geq&\liminf_{t\to\infty}\min_{x\neq x^*(y)}-\frac{1}{t^{\frac{1}{1+\xi}}}\ln[2\\
	&\max\{\exp(\eta^{-\xi}(-\eta ^{1+\xi}\ln\eta ^{1+\xi}+\eta ^{1+\xi}- 1 ) (C\lambda_L)^\frac{1}{1+\xi}t^\frac{1}{1+\xi}+\ln 2),\exp(-t^{\frac{1}{1+\xi}} H_3^{(t)}(x,y;\eta))\}]\\
	=&\liminf_{t\to\infty}\min_{x\neq x^*(y)}\min\{ \eta^{-\xi}(-\eta ^{1+\xi}\ln\eta ^{1+\xi}+\eta ^{1+\xi}- 1 ) (C\lambda_L)^\frac{1}{1+\xi}+\frac{\ln 4}{t^\frac{1}{1+\xi}},         H_3^{(t)}(x,y;\eta)+\frac{\ln 2}{t^\frac{1}{1+\xi}}   \}\\
	=&\min_{x\neq x^*(y)}\min\{\eta^{-\xi}(\eta ^{1+\xi}\ln\eta ^{1+\xi}-\eta ^{1+\xi}+ 1 ) (C\lambda_L)^\frac{1}{1+\xi},H_3(x,y;\eta)\}\\
	=& (C \lambda_L)^\frac{1}{1+\xi}\min_{x\neq x^*(y)}\min\{\eta^{-\xi}(\eta ^{1+\xi}\ln\eta ^{1+\xi}-\eta ^{1+\xi}+ 1 ),\eta \frac{(f_x(y)-f_{x^*}(y))^2}{2({\sigma_{x^*}^2(y)}+{\sigma_{x}^2(y)})}\}
\end{split}
\end{equation}
The first inequality holds due to (\ref{equation:CDLupper})	. The first equality holds as $\max_{x\neq x^*(y)}d({f^\dagger}^{(t)}(y),{f^*(y)})$ is a globally bounded value and any bounded value will vanish under operations of $-\frac{1}{t^\frac{1}{1+\xi}}\ln[\cdot]$. The second equality holds by changing of the signature. The second inequality holds due to (\ref{inequality: mutual probability}). The third inequality holds due to monotonicity of $\ln$ function. The fourth equality holds due to each of the finite number of terms in $\min$ converges to a constant and thus $\liminf$ and $\min$ is interchangeable.  Notice that the inequality holds for arbitrary $\eta\in(0,1)$, $\forall x\in\chi, y\in Y^\circ$. Therefore, to achieve the tightest bound:
\begin{equation}
\begin{split}
	&\liminf_{t\to\infty}-\frac{1}{t^{\frac{1}{1+\xi}}}\ln \upsilon^{(t)} (y)\\
	\geq&\max_{\eta\in(0,1),\forall x\in\chi}\big\{(C \lambda_L)^\frac{1}{1+\xi}\min_{x\neq x^*(y)}\min\{\eta^{-\xi}(\eta ^{1+\xi}\ln\eta ^{1+\xi}-\eta ^{1+\xi}+ 1 ),\eta \frac{(f_x(y)-f_{x^*}(y))^2}{2({\sigma_{x^*}^2(y)}+{\sigma_{x}^2(y)})}\}\big\}\\
	=&(C \lambda_L)^\frac{1}{1+\xi}\min_{x\neq x^*(y)} \frac{\left(f_x(y)-f_{x^*}(y)\right)^2}{2 \left(\sigma _{x^*}^2(y)+\sigma _x^2(y)\right)} \left(-\frac{1}{W_{-1}\left(-\exp \left(-\frac{\left(f_x(y)-f_{x^*}(y)\right)^2}{2
			\left(\sigma _{x^*}^2(y)+\sigma _x^2(y)\right)}-1\right)\right)}\right)^{\frac{1}{1+\xi}}
\end{split}\label{inequality: Lambert W function}
\end{equation}
where $W_{-1}$ is the lower branch of Lambert $W$ function (i.e. $w=W_{-1}(z)$ if $we^w=z$ for $z\in[-e^{-1},0)]$). The unique optimal $\eta$ is achieved inside the domain as it is easy to show that $\eta^{-\xi}(\eta ^{1+\xi}\ln\eta ^{1+\xi}-\eta ^{1+\xi}+ 1 )$ is a strictly monotone decreasing function w.r.t $\eta$ and $\eta \frac{(f_x(y)-f_{x^*}(y))^2}{2({\sigma_{x^*}^2(y)}+{\sigma_{x}^2(y)})}$ is strictly monotone increasing function w.r.t $\eta$. Two functions will have a unique intersection for an $\eta\in(0,1)$. Define $\zeta$ function for $y\in Y^\circ$:
\begin{equation}
\zeta(y):= \min_{x\neq x^*(y)}\frac{\left(f_x(y)-f_{x^*}(y)\right)^2}{2 \left(\sigma _{x^*}^2(y)+\sigma _x^2(y)\right)},
\end{equation}
Notice that the following function is monotone increasing for $z\in(0,\infty)$:
\begin{equation}
z \left(-\frac{1}{W_{-1}\left(-\exp \left(-z-1\right)\right)}\right)^{\frac{1}{1+\xi}}
\end{equation}
Thus, based on (\ref{inequality: Lambert W function}):
\begin{equation}
\begin{split}
	&\liminf_{t\to\infty}-\frac{1}{t^{\frac{1}{1+\xi}}}\ln \upsilon^{(t)} (y)\\
	\geq&(C \lambda_L)^\frac{1}{1+\xi}\zeta(y) \left(-\frac{1}{W_{-1}\left(-\exp \left(-\zeta(y)-1\right)\right)}\right)^{\frac{1}{1+\xi}}\\
	=:&{H'_y}
\end{split}\label{inequality: rate of CDL}
\end{equation}
\endproof{$\blacksquare$}

\subsection{Proposition \ref{proposition: subset rate}}
\proof{Proof}
\begin{align*}
&\liminf_{t\to\infty}-\frac{1}{t^{\frac{1}{1+\xi}}}\ln[\int_{y\in Y_{sub}}\upsilon^{(t)} (y)dP_Y(y)]\\
\geq&\liminf_{t\to\infty}-\frac{1}{t^{\frac{1}{1+\xi}}}\ln[\int_{y\in Y_{sub}}dP_Y(y)\max_{y\in Y_{sub}}\upsilon^{(t)} (y)]\\
=&\liminf_{t\to\infty}-\frac{1}{t^{\frac{1}{1+\xi}}}\ln[\max_{y\in Y_{sub}}\upsilon ^{(t)}(y)]\\
=&\liminf_{t\to\infty}\min_{y\in Y_{sub}}-\frac{1}{t^{\frac{1}{1+\xi}}}\ln[\upsilon ^{(t)}(y)]\\		
=&\min_{y\in Y_{sub}}\liminf_{t\to\infty}-\frac{1}{t^{\frac{1}{1+\xi}}}\ln[\upsilon ^{(t)}(y)]\\		
\geq&\min_{y\in Y_{sub}} {H'_y}\\
=:&H_{Y_{sub}}, 
\end{align*}
The first inequality holds as $Y_{sub}$ is a compact set, indicating that the maximum can be attained. The third equality holds as $Y_{sub}$ is a compact set and thus $\liminf_{t\to\infty}$ and $\min_{y\in Y_{sub}}$ can be interchanged. The second inequality holds according to (\ref{inequality: rate of CDL})
\endproof{$\blacksquare$}

\subsection{Lemma \ref{Lemma:mt static}}
\proof{Proof}
Since the discussion is for each $x, y$, we omit specification of $x, y$ if no confusion is caused. We prove by using Borel-Cantelli Lemma, which is to show that $\forall \eta\in(0,1)$, 
\begin{align*}
\sum_{N=1}^\infty P\{\Psi^{(t)} <\lfloor\eta\psi^{(t)} \rfloor\}<\infty.
\end{align*}
Notice that by definition:
\begin{equation}
\psi^{(t)}:=\max\{m\in[1,t],m\in\mathbb{N}:(Ca_{\inf, m})^\frac{1}{1+\xi}t^\frac{1}{1+\xi}\geq m\},\label{definition: psi}
\end{equation}
where \begin{equation}
\begin{split}
	C:=&\frac{\pi^\frac{d}{2}}{\Gamma(\frac{d}{2}+1)}c^d\\
	a_{\inf, m}:=&\inf_{y'\in B(r_m;y)} a(x,y').
\end{split}
\end{equation}
The definition of $	\psi^{(t)}$ immediately suggests the following inequality:
\begin{equation}
(Ca_{\inf, \psi^{(t)}+1})^\frac{1}{1+\xi}t^\frac{1}{1+\xi}-1
<{\psi^{(t)}}
\leq (Ca_{\inf, \psi^{(t)}})^\frac{1}{1+\xi}t^\frac{1}{1+\xi}\label{inequality:  psi}
\end{equation}
The first inequality is due to the fact that if $\psi^{(t)}=m$, $m+1$ is not qualified to be $\psi^{(t)}$. The second inequality is derived based on the definition of $\psi^{(t)}$. Notice that for a fixed $t$, $a_{\inf, m}$ as a sequence of $m$ is monotone increasing. This monotonicity suggests that: 
\begin{equation}
{\psi^{(t)}}>(Ca_{\inf, \psi^{(t)}+1})^\frac{1}{1+\xi}t^\frac{1}{1+\xi}-1\geq (Ca_{\inf, 1})^\frac{1}{1+\xi}t^\frac{1}{1+\xi}-1.
\end{equation}
Since $(Ca_{\inf, 1})^\frac{1}{1+\xi}t^\frac{1}{1+\xi}-1$ goes to infinity as $t\to\infty$, we have that as $t\to\infty$, $${\psi^{(t)}} \to\infty$$ and $$r_{{\psi^{(t)}} }=c{\psi^{(t)}}^{-\xi}\to 0.$$ Therefore, for $y\in Y^\circ$, $\exists T^*$ such that $\forall t\geq T^*$, $B(r_{{\psi^{(t)}} })$ is a ball inside $Y^\circ$. In the following discussion, we assume $t\geq T^*$ so that $B(r_{{\psi^{(t)}} })\subset Y^\circ$. 

Since $\lim_{t\to\infty}{\psi^{(t)}} =\infty$, as $t\to\infty$, formulation (\ref{inequality:  psi}) leads to:
\begin{equation}
(Ca)^\frac{1}{1+\xi}	\leftarrow(Ca_{\inf, \psi^{(t)}+1})^\frac{1}{1+\xi}-\frac{1}{t^\frac{1}{1+\xi}}
<\frac{\psi^{(t)}}{t^\frac{1}{1+\xi}}
\leq (Ca_{\inf, \psi^{(t)}})^\frac{1}{1+\xi}\to(Ca)^\frac{1}{1+\xi}.\label{inequality:  lim}
\end{equation}
The convergence holds due to $a(x,y)$ is continuous on $Y^\circ$ for any fixed $x$, leading to $\lim_{t\to\infty}a_{\inf,x,y, \psi^{(t)}}=a(x,y)$. By squeeze theorem,
\begin{equation}
\lim_{t\to\infty}\frac{\psi^{(t)}}{t^\frac{1}{1+\xi}}=(Ca)^\frac{1}{1+\xi} \label{limit:ratio}
\end{equation}

On stage $t$, denote the probability such that a point $(x',y')$ randomly generated by $a(x,y)$ falls into the neighborhood of $(x,y)$ with radius of $r_m$ (e.g. $x'=x, y'\in B(r_m;y)$) as $p_m$ (again, $x,y$ is omitted in this notation).  Notice that in the following discussion, $m$ is allowed to take a continuous value and $r_m=cm^{-\xi/d}$ still applies. Define:
\begin{equation}
a_{\sup,m}:=\sup_{y'\in B(r_m;y)} a(x,y')
\end{equation}For $t\geq T^*$, we have the following bound for $p_{m}$ for arbitrary $m>0$:
\begin{equation}
Cm^{-\xi} a_{\sup,m}=\frac{\pi^\frac{d}{2}}{\Gamma(\frac{d}{2}+1)}r_{m}^d a_{\sup,m}\geq p_{m}\geq\frac{\pi^\frac{d}{2}}{\Gamma(\frac{d}{2}+1)}r_{m}^d a_{\inf,m}=Cm^{-\xi} a_{\inf,m}\label{inequality: p bound}
\end{equation}
The inequality holds as $p_{m}$ is the integral of $a(x,y)$ on $B(r_{m};y)$ and is thus bounded by integral of $a_{\sup,m}$ and that of $a_{\inf,m}$.

Let $n_{m}^{(t)}$ (i.e. $n_{m}^{(t)}(x,y)$) denotes the total number of samples fall into $B (r_{m};y)$ at time $t$. For $u<0$, we have the following upper bound:
\begin{align*}
&P\{\Psi ^{(t)}<\lfloor\eta{\psi^{(t)}} \rfloor\}\\
=&P\{ n_{\lfloor\eta{\psi^{(t)}} \rfloor}^{(t)}<\lfloor\eta{\psi^{(t)}} \rfloor\}\\
\leq	&P\{ n_{{\eta\psi^{(t)}}}^{(t)}<\eta{\psi^{(t)}}  \}\\
\leq	&\inf_{u<0}(1- p_{\eta \psi^{(t)}}+ p_{\eta \psi^{(t)}}\exp(u))^t\exp(-u\eta \psi^{(t)})\\
=&\exp(-\eta \psi^{(t)}\ln(\frac{\eta \psi^{(t)}}{ t  p_{\eta \psi^{(t)}} })-( t -\eta \psi^{(t)})\ln(\frac{t-\eta \psi^{(t)}}{ t-t  p_{\eta \psi^{(t)}} })  )\\
\leq&\exp(-\eta \psi^{(t)}\ln(\frac{\eta \psi^{(t)}}{ t  C(\eta\psi^{(t)})^{-\xi} a_{\inf,\eta\psi^{(t)}} })-( t -\eta \psi^{(t)})\ln(\frac{t-\eta \psi^{(t)}}{ t-t  C(\eta\psi^{(t)})^{-\xi} a_{\inf,\eta\psi^{(t)}}})  )
\end{align*}
This leads to the fact that:
\begin{equation}
\begin{split}
	&\lim_{t\to\infty}-\frac{1}{t^\frac{1}{1+\xi}}\ln(P\{\Psi ^{(t)}<\lfloor\eta{\psi^{(t)}} \rfloor\})\\
	\geq&\lim_{t\to\infty}-\frac{1}{t^\frac{1}{1+\xi}}\ln(\exp(-\eta \psi^{(t)}\ln(\frac{\eta \psi^{(t)}}{ t  C(\eta\psi^{(t)})^{-\xi} a_{\inf,\eta\psi^{(t)}} })-( t  -\eta \psi^{(t)})\ln(\frac{t-\eta \psi^{(t)}}{ t-t  C(\eta\psi^{(t)})^{-\xi} a_{\inf,\eta\psi^{(t)}}})  ))\\
	=&\lim_{t\to\infty}\frac{1}{t^\frac{1}{1+\xi}}(\eta \psi^{(t)}\ln(\frac{\eta \psi^{(t)}}{ t  C(\eta\psi^{(t)})^{-\xi} a_{\inf,\eta\psi^{(t)}} })+( t  -\eta \psi^{(t)})\ln(\frac{t-\eta \psi^{(t)}}{ t-t  C(\eta\psi^{(t)})^{-\xi} a_{\inf,\eta\psi^{(t)}}})  )\\
	=&\eta (Ca)^\frac{1}{1+\xi}\ln(\eta^{1+\xi})+\lim_{t\to\infty}\frac{ t  -\eta \psi^{(t)}}{t^\frac{1}{1+\xi}}\ln(\frac{t-\eta \psi^{(t)}}{ t-t  C(\eta\psi^{(t)})^{-\xi} a_{\inf,\eta\psi^{(t)}}})  \\
	=&\eta (Ca)^\frac{1}{1+\xi}\ln(\eta^{1+\xi})+\lim_{t\to\infty}\frac{ t  -\eta \psi^{(t)}}{t^\frac{1}{1+\xi}}\frac{t  C(\eta\psi^{(t)})^{-\xi} a_{\inf,\eta\psi^{(t)}}-\eta \psi^{(t)}}{ t-t  C(\eta\psi^{(t)})^{-\xi} a_{\inf,\eta\psi^{(t)}}}  \\
	=&\eta (Ca)^\frac{1}{1+\xi}\ln(\eta^{1+\xi})+\lim_{t\to\infty}(\eta^{-\xi}-\eta)(Ca)^\frac{1}{1+\xi}\frac{t  -\eta \psi^{(t)}}{ t-t  C(\eta\psi^{(t)})^{-\xi} a_{\inf,\eta\psi^{(t)}}} \\
	=&\eta (Ca)^\frac{1}{1+\xi}(\ln\eta^{1+\xi}+\eta^{-1-\xi}-1)
\end{split}\label{inequality: insufficient sampling lower bound}
\end{equation}
On the other hand:
\begin{align*}
&P\{\Psi ^{(t)}<\lfloor\eta{\psi^{(t)}} \rfloor\}\\
=&P\{ n_{\lfloor\eta{\psi^{(t)}} \rfloor}^{(t)}<\lfloor\eta{\psi^{(t)}} \rfloor\}\\
\geq	&P\{ n_{\lfloor\eta{\psi^{(t)}} \rfloor-1}^{(t)}<\lfloor\eta{\psi^{(t)}} \rfloor-1 \}\\
\geq&\frac{1}{\sqrt{8t  \frac{\lfloor\eta{\psi^{(t)}} \rfloor-1}{t}(1-\frac{\lfloor\eta{\psi^{(t)}} \rfloor-1}{t})}}\\
&\exp(-(\lfloor\eta{\psi^{(t)}} \rfloor-1)\ln(\frac{\lfloor\eta{\psi^{(t)}} \rfloor-1}{ t  p_{\lfloor\eta{\psi^{(t)}} \rfloor-1} })-( t-(\lfloor\eta{\psi^{(t)}} \rfloor-1))\ln(\frac{t-(\lfloor\eta{\psi^{(t)}} \rfloor-1)}{ t-t  p_{\lfloor\eta{\psi^{(t)}} \rfloor-1} })  )\\
=&\frac{1}{\sqrt{8t  \frac{\lfloor\eta{\psi^{(t)}} \rfloor-1}{t}(1-\frac{\lfloor\eta{\psi^{(t)}} \rfloor-1}{t})}}\\
&\exp[-t^{\frac{1}{1+\xi}}(\frac{\lfloor\eta{\psi^{(t)}} \rfloor-1}{t^{\frac{1}{1+\xi}}}\ln(\frac{\frac{\lfloor\eta{\psi^{(t)}} \rfloor-1}{t^{\frac{1}{1+\xi}}}}{ t^\frac{\xi}{1+\xi}  	C(\lfloor\eta{\psi^{(t)}} \rfloor-1)^{-\xi} a_{\sup,\lfloor\eta{\psi^{(t)}} \rfloor-1} })\\
&+( t^\frac{\xi}{1+\xi}-\frac{\lfloor\eta{\psi^{(t)}} \rfloor-1}{t^{\frac{1}{1+\xi}}})\ln(1+\frac{t  C(\lfloor\eta{\psi^{(t)}} \rfloor-1)^{-\xi} a_{\sup,\lfloor\eta{\psi^{(t)}} \rfloor-1}-(\lfloor\eta{\psi^{(t)}} \rfloor-1)}{ t-t  C(\lfloor\eta{\psi^{(t)}} \rfloor-1)^{-\xi} a_{\sup,\lfloor\eta{\psi^{(t)}} \rfloor-1} }))]
\end{align*}
This leads to 
\begin{align*}
&\lim_{t\to\infty}-\frac{1}{t^\frac{1}{1+\xi}}\ln(P\{\Psi ^{(t)}<\lfloor\eta{\psi^{(t)}} \rfloor\})\\
\leq&\lim_{t\to\infty}-\frac{1}{t^\frac{1}{1+\xi}}\ln\Big[\frac{1}{\sqrt{8t  \frac{\lfloor\eta{\psi^{(t)}} \rfloor-1}{t}(1-\frac{\lfloor\eta{\psi^{(t)}} \rfloor-1}{t})}}\\
&\exp[-t^{\frac{1}{1+\xi}}(\frac{\lfloor\eta{\psi^{(t)}} \rfloor-1}{t^{\frac{1}{1+\xi}}}\ln(\frac{\frac{\lfloor\eta{\psi^{(t)}} \rfloor-1}{t^{\frac{1}{1+\xi}}}}{ t^\frac{\xi}{1+\xi}  	C(\lfloor\eta{\psi^{(t)}} \rfloor-1)^{-\xi} a_{\sup,\lfloor\eta{\psi^{(t)}} \rfloor-1} })\\
&+( t^\frac{\xi}{1+\xi}-\frac{\lfloor\eta{\psi^{(t)}} \rfloor-1}{t^{\frac{1}{1+\xi}}})\ln(1+\frac{t  C(\lfloor\eta{\psi^{(t)}} \rfloor-1)^{-\xi} a_{\sup,\lfloor\eta{\psi^{(t)}} \rfloor-1}-(\lfloor\eta{\psi^{(t)}} \rfloor-1)}{ t-t  C(\lfloor\eta{\psi^{(t)}} \rfloor-1)^{-\xi} a_{\sup,\lfloor\eta{\psi^{(t)}} \rfloor-1} }))]\Big]\\
=&\lim_{t\to\infty}-\frac{1}{t^\frac{1}{1+\xi}}  \ln\Big[\exp[-t^{\frac{1}{1+\xi}}(\frac{\lfloor\eta{\psi^{(t)}} \rfloor-1}{t^{\frac{1}{1+\xi}}}\ln(\frac{\frac{\lfloor\eta{\psi^{(t)}} \rfloor-1}{t^{\frac{1}{1+\xi}}}}{ t^\frac{\xi}{1+\xi}  	C(\lfloor\eta{\psi^{(t)}} \rfloor-1)^{-\xi} a_{\sup,\lfloor\eta{\psi^{(t)}} \rfloor-1} })\\
&+( t^\frac{\xi}{1+\xi}-\frac{\lfloor\eta{\psi^{(t)}} \rfloor-1}{t^{\frac{1}{1+\xi}}})\ln(1+\frac{t  C(\lfloor\eta{\psi^{(t)}} \rfloor-1)^{-\xi} a_{\sup,\lfloor\eta{\psi^{(t)}} \rfloor-1}-(\lfloor\eta{\psi^{(t)}} \rfloor-1)}{ t-t  C(\lfloor\eta{\psi^{(t)}} \rfloor-1)^{-\xi} a_{\sup,\lfloor\eta{\psi^{(t)}} \rfloor-1} }))]\Big]\\
=&\lim_{t\to\infty}\frac{\lfloor\eta{\psi^{(t)}} \rfloor-1}{t^{\frac{1}{1+\xi}}}\ln(\frac{\frac{\lfloor\eta{\psi^{(t)}} \rfloor-1}{t^{\frac{1}{1+\xi}}}}{ t^\frac{\xi}{1+\xi}  	C(\lfloor\eta{\psi^{(t)}} \rfloor-1)^{-\xi} a_{\sup,\lfloor\eta{\psi^{(t)}} \rfloor-1} })\\
&+( t^\frac{\xi}{1+\xi}-\frac{\lfloor\eta{\psi^{(t)}} \rfloor-1}{t^{\frac{1}{1+\xi}}})\ln(1+\frac{t  C(\lfloor\eta{\psi^{(t)}} \rfloor-1)^{-\xi} a_{\sup,\lfloor\eta{\psi^{(t)}} \rfloor-1}-(\lfloor\eta{\psi^{(t)}} \rfloor-1)}{ t-t  C(\lfloor\eta{\psi^{(t)}} \rfloor-1)^{-\xi} a_{\sup,\lfloor\eta{\psi^{(t)}} \rfloor-1} })\\
=&\eta (Ca)^\frac{1}{1+\xi}\ln \eta^{1+\xi}+\lim_{t\to\infty}( t^\frac{\xi}{1+\xi}-\frac{\lfloor\eta{\psi^{(t)}} \rfloor-1}{t^{\frac{1}{1+\xi}}})\frac{t  C(\lfloor\eta{\psi^{(t)}} \rfloor-1)^{-\xi} a_{\sup,\lfloor\eta{\psi^{(t)}} \rfloor-1}-(\lfloor\eta{\psi^{(t)}} \rfloor-1)}{ t-t  C(\lfloor\eta{\psi^{(t)}} \rfloor-1)^{-\xi} a_{\sup,\lfloor\eta{\psi^{(t)}} \rfloor-1} }\\
=&\eta (Ca)^\frac{1}{1+\xi}(\ln\eta^{1+\xi}+\eta^{-1-\xi}-1)
\end{align*}
Based on the above and (\ref{inequality: insufficient sampling lower bound}):
\begin{equation}
\lim_{t\to\infty}-\frac{1}{t^\frac{1}{1+\xi}}\ln(P\{\Psi ^{(t)}<\lfloor\eta{\psi^{(t)}} \rfloor\})=\eta (Ca)^\frac{1}{1+\xi}(\ln\eta^{1+\xi}+\eta^{-1-\xi}-1)
\end{equation}
Therefore, for arbitrary $\delta>0$ small enough such that $\eta (Ca)^\frac{1}{1+\xi}(\ln\eta^{1+\xi}+\eta^{-1-\xi}-1)-\delta>0$, exists $t>T^{**}$ such that:
\begin{equation}
P\{\Psi ^{(t)}<\lfloor\eta{\psi^{(t)}} \rfloor\}\leq \exp[-[\eta (Ca)^\frac{1}{1+\xi}(\ln\eta^{1+\xi}+\eta^{-1-\xi}-1)-\delta](t^\frac{1}{1+\xi})]
\end{equation}
Therefore:
\begin{equation}\\
\begin{split}
	&\lim_{t\to\infty}\sum_{t'=\max\{T^*,T^{**}\}}^tP\{\Psi <\lfloor{\eta\psi^{(t')}} \rfloor\}\\
	\leq&\lim_{t\to\infty}\sum_{t'=\max\{T^*,T^{**}\}}^{t}\exp[-[\eta (Ca)^\frac{1}{1+\xi}(\ln\eta^{1+\xi}+\eta^{-1-\xi}-1)-\delta]({t'}^\frac{1}{1+\xi})]<\infty
\end{split}
\end{equation}
Based on Borel-Cantelli lemma, $\Psi ^{(t)}<\lfloor\eta{\psi^{(t)}} \rfloor$ holds for finitely many $t$ with probability $1$. Thus $\Psi ^{(t)}\geq\lfloor\eta{\psi^{(t)}} \rfloor$ holds for all but finitely many $t$ with probability $1$.
\endproof{$\blacksquare$}

\subsection{Theorem \ref{theorem: CDL convergence rate static}}
\proof{Proof}
$\forall y\in Y\backslash MO$, $\upsilon^{(t)} $ is upper bounded by:
\begin{equation}
\begin{split}
	\upsilon^{(t)} (y)=&E_{\mathcal{D}^{(t)} }d({f^\dagger}^{(t)}(y),{f^*(y)})\\
	=&\sum_{x\neq x^*(y)}d({f^\dagger}^{(t)}(y),{f^*(y)})P\{\hat{x}^{*(t)}(y)=x\}\\
	\leq&(|\chi|-1)\max_{x\neq x^*(y)}d({f^\dagger}^{(t)}(y),{f^*(y)})P\{{\hat{f}_x^{(t)}(y)}<\hat{f}_{x^*}^{(t)}(y)\}
\end{split},\label{equation:CDLupper static}
\end{equation}
Now we focus on $P\{{\hat{f}_x^{(t)}(y)}<\hat{f}_{x^*}^{(t)}(y)\}$. For each $x\neq x^*(y)$, this probability can be separated into the one caused by observation noise and the one caused by bias. For $x\neq x^*(y)$:
\begin{align*}
&P\{\hat{f}_{x}^{(t)}(y)<\hat{f}_{x^*}^{(t)}(y)\}\\
=&P\{\bar{\epsilon}^{(t)} _{x, y}+ \frac{1}{{\Psi^{(t)}_{x,y}}}\sum_{i:x_{[i]}=x, y_{[i]}\in {b _{x, y}}}[f_{[i]}-f_{x}(y)]+f_{x}(y)   \\
&<{\bar{\epsilon}^{(t)} _{x^*, y}}+ \frac{1}{{\Psi_{x^*, y}^{(t)}}}\sum_{i:x_{[i]}= x^*(y), y_{[i]}=b _{x^*, y}}[f_{[i]}-{f_{x^*}(y)}]+{f_{x^*}(y)}  \}\\
\leq &P\{{\bar{\epsilon}^{(t)} _{x^*, y}}-\bar{\epsilon}^{(t)} _{x, y}> -({r_{\Psi_{x^*, y}^{(t)}}}+{r_{\Psi_{x, y}^{(t)}}})L_1  + f_{x}(y)-{f_{x^*}(y)} \}  \\
=&P\{{\bar{\epsilon}^{(t)} _{x^*, y}}-\bar{\epsilon}^{(t)} _{x, y}> -({r_{\Psi_{x^*, y}^{(t)}}}+{r_{\Psi_{x, y}^{(t)}}})L_1  + f_{x}(y)-{f_{x^*}(y)},\{ \Psi_{x^*, y}^{(t)}< \lfloor\eta_{x^*, y}\psi_{x^*, y} ^{(t)}\rfloor\}\cup\{\Psi_{x, y}^{(t)}<  \lfloor\eta_{x, y}\psi_{x, y} ^{(t)}\rfloor\} \} \\
&+P\{{\bar{\epsilon}^{(t)} _{x^*, y}}-\bar{\epsilon}^{(t)} _{x, y}> -({r_{\Psi_{x^*, y}^{(t)}}}+{r_{\Psi_{x, y}^{(t)}}})L_1  + f_{x}(y)-{f_{x^*}(y)},\{ \Psi_{x^*, y}^{(t)}\geq \lfloor\eta_{x^*, y}\psi_{x^*, y} ^{(t)}\rfloor\}\cap\{\Psi_{x, y}^{(t)}\geq  \lfloor\eta_{x, y}\psi_{x, y} ^{(t)}\rfloor\} \},
\end{align*}
where $\eta_{x, y}\in(0,1)$ is now associated with each $(x,y)$. For the first term, it is upper bounded by:
\begin{equation}
\begin{split}
	&P\{{\bar{\epsilon}^{(t)} _{x^*, y}}-\bar{\epsilon}^{(t)} _{x, y}> -({r_{\Psi_{x^*, y}^{(t)}}}+{r_{\Psi_{x, y}^{(t)}}})L_1  + f_{x}(y)-{f_{x^*}(y)},\{ \Psi_{x^*, y}^{(t)}< \lfloor\eta_{x^*, y}\psi_{x^*, y} ^{(t)}\rfloor\}\cup\{\Psi_{x, y}^{(t)}<  \lfloor\eta_{x, y}\psi_{x, y} ^{(t)}\rfloor\} \}\\
	\leq &P\{ \Psi_{x^*, y}^{(t)}< \lfloor\eta_{x^*, y}\psi_{x^*, y} ^{(t)}\rfloor\}+P\{ \Psi_{x, y}^{(t)}<  \lfloor\eta_{x, y}\psi_{x, y} ^{(t)}\rfloor\}.
\end{split}\label{inequality: insufficient samples}\end{equation}
For the second term:
\begin{equation}
\begin{split}
	&P\{{\bar{\epsilon}^{(t)} _{x^*, y}}-\bar{\epsilon}^{(t)} _{x, y}> -({r_{\Psi_{x^*, y}^{(t)}}}+{r_{\Psi_{x, y}^{(t)}}})L_1  + f_{x}(y)-{f_{x^*}(y)},\{ \Psi_{x^*, y}^{(t)}\geq \lfloor\eta_{x^*, y}\psi_{x^*, y} ^{(t)}\rfloor\}\cap\{\Psi_{x, y}^{(t)}\geq  \lfloor\eta_{x, y}\psi_{x, y} ^{(t)}\rfloor\} \} \\
	\leq&\sup_{\Psi_{x^*, y}^{(t)}\geq \lfloor\eta_{x^*, y}\psi_{x^*, y} ^{(t)}\rfloor, \Psi_{x, y}^{(t)}\geq \lfloor\eta_{x, y}\psi_{x, y} ^{(t)}\rfloor,(\sigma_{[i_1]},\sigma_{[i_2]},\dots,\sigma_{[i_{\Psi_{x^*, y}}]}),(\sigma_{[k_1]},\sigma_{[k_2]},\dots,\sigma_{[k_{\Psi^{(t)}_{x,y}}]})}\\
	&P\{{\bar{\epsilon}^{(t)} _{x^*, y}}-\bar{\epsilon}^{(t)} _{x, y}> -({r_{\Psi_{x^*, y}^{(t)}}}+{r_{\Psi_{x, y}^{(t)}}})L_1  + f_{x}(y)-{f_{x^*}(y)} \Big|(\sigma_{[i_1]},\sigma_{[i_2]},\dots,\sigma_{[i_{\Psi_{x^*, y}^{(t)}}]}),(\sigma_{[k_1]},\sigma_{[k_2]},\dots,\sigma_{[k_{\Psi^{(t)}_{x,y}}]}),\\
	&{\Psi_{x^*, y}^{(t)}}, {\Psi_{x, y}^{(t)} }\}\\
	=&\sup_{\Psi_{x^*, y}^{(t)}\geq \lfloor\eta_{x^*, y}\psi_{x^*, y} ^{(t)}\rfloor, \Psi_{x, y}^{(t)}\geq \lfloor\eta_{x, y}\psi_{x, y} ^{(t)}\rfloor,(\sigma_{[i_1]},\sigma_{[i_2]},\dots,\sigma_{[i_{\Psi_{x^*, y}}]}),(\sigma_{[k_1]},\sigma_{[k_2]},\dots,\sigma_{[k_{\Psi^{(t)}_{x,y}}]})}\\
	&\Phi(-\frac{-({r_{\Psi_{x^*, y}^{(t)}}}+{r_{\Psi_{x, y}^{(t)}}})L_1  + f_{x}(y)-{f_{x^*}(y)}}{\sqrt{\frac{\sum_{j=1}^{{\Psi_{x^*,y}}^{(t)}}\sigma_{x_{[i_j]}}^2(y_{[i_j]})}{({\Psi_{x^*,y}}^{(t)})^2}+\frac{\sum_{j=1}^{{\Psi_{x,y}}^{(t)}}\sigma_{x_{[i_k]}}^2(y_{[i_k]})}{({\Psi_{x,y}}^{(t)})^2}}})\\
	=&\Phi(-\frac{-({r_{ \lfloor\eta_{x^*, y}\psi_{x^*, y}^{(t)}\rfloor}}+{r_{ \lfloor\eta_{x, y}\psi_{x, y}^{(t)}\rfloor}})L_1  + f_{x}(y)-{f_{x^*}(y)}}{\sqrt{\frac{(\sigma_{x^*}(y)+L_2c\lfloor\eta_{x^*, y}\psi_{x^*, y}^{(t)}\rfloor^{-\frac{\xi}{d}})^2}
			{\lfloor\eta_{x^* y}\psi_{x^*, y}^{(t)}\rfloor}
			+\frac{(\sigma_{x}(y)+L_2c\lfloor\eta_{x, y}\psi_{x, y}^{(t)}\rfloor^{-\frac{\xi}{d}})^2}
			{\lfloor\eta_{x, y}\psi_{x, y}^{(t)}\rfloor}}})\\
	\leq&\exp(-\frac{\left(-({r_{ \lfloor\eta_{x^*, y}\psi_{x^*, y}^{(t)}\rfloor}}+{r_{ \lfloor\eta_{x, y}\psi_{x, y}^{(t)}\rfloor}})L_1  + f_{x}(y)-{f_{x^*}(y)}\right)^2}{2\left({\frac{(\sigma_{x^*}(y)+L_2c\lfloor\eta_{x^*, y}\psi_{x^*, y}^{(t)}\rfloor^{-\frac{\xi}{d}})^2}
			{\lfloor\eta_{x^* y}\psi_{x^*, y}^{(t)}\rfloor}
			+\frac{(\sigma_{x}(y)+L_2c\lfloor\eta_{x, y}\psi_{x, y}^{(t)}\rfloor^{-\frac{\xi}{d}})^2}
			{\lfloor\eta_{x, y}\psi_{x, y}^{(t)}\rfloor}}\right)})
\end{split}\label{inequation: mixed probability upper bound static}
\end{equation}

Therefore, we have:
\begin{equation}
\begin{split}
	&\liminf_{t\to\infty}-\frac{1}{t^{\frac{1}{1+\xi}}}\ln \upsilon^{(t)} (y)\\
	\geq&\liminf_{t\to\infty}-\frac{1}{t^{\frac{1}{1+\xi}}}\ln\Big[(|\chi|-1)\max_{x\neq x^*(y)}P\{{\hat{f}_x^{(t)}(y)}<\hat{f}_{x^*}^{(t)}(y)\}\max_{x\neq x^*(y)}d({f^\dagger}^{(t)}(y),{f^*(y)})\Big]\\
	=&\liminf_{t\to\infty}-\frac{1}{t^{\frac{1}{1+\xi}}}\ln \max_{x\neq x^*(y)} P\{\hat{f}_x^{(t)}(y)<\hat{f}_{x^*}^{(t)}(y)\}\\
	=&\liminf_{t\to\infty}\min_{x\neq x^*(y)}-\frac{1}{t^{\frac{1}{1+\xi}}}\ln  P\{\hat{f}_x^{(t)}(y)<\hat{f}_{x^*}^{(t)}(y)\}\\
	\geq&\liminf_{t\to\infty}\min_{x\neq x^*(y)}-\frac{1}{t^{\frac{1}{1+\xi}}}\ln \Big[P\{ \Psi_{x^*, y}^{(t)}< \lfloor\eta_{x^*, y}\psi_{x^*, y} ^{(t)}\rfloor\}+P\{ \Psi_{x, y}^{(t)}<  \lfloor\eta_{x, y}\psi_{x, y} ^{(t)}\rfloor\}\\
	&+\exp(-\frac{\left(-({r_{ \lfloor\eta_{x^*, y}\psi_{x^*, y}^{(t)}\rfloor}}+{r_{ \lfloor\eta_{x, y}\psi_{x, y}^{(t)}\rfloor}})L_1  + f_{x}(y)-{f_{x^*}(y)}\right)^2}{2\left({\frac{(\sigma_{x^*}(y)+L_2c\lfloor\eta_{x^*, y}\psi_{x^*, y}^{(t)}\rfloor^{-\frac{\xi}{d}})^2}
			{\lfloor\eta_{x^* y}\psi_{x^*, y}^{(t)}\rfloor}
			+\frac{(\sigma_{x}(y)+L_2c\lfloor\eta_{x, y}\psi_{x, y}^{(t)}\rfloor^{-\frac{\xi}{d}})^2}
			{\lfloor\eta_{x, y}\psi_{x, y}^{(t)}\rfloor}}\right)})\Big]\\
	\geq&\liminf_{t\to\infty}\min_{x\neq x^*(y)}-\frac{1}{t^{\frac{1}{1+\xi}}}\ln \Big[3\max\{P\{ \Psi_{x^*, y}^{(t)}< \lfloor\eta_{x^*, y}\psi_{x^*, y} ^{(t)}\rfloor\}, P\{ \Psi_{x, y}^{(t)}<  \lfloor\eta_{x, y}\psi_{x, y} ^{(t)}\rfloor\},\\
	&\exp(-\frac{\left(-({r_{ \lfloor\eta_{x^*, y}\psi_{x^*, y}^{(t)}\rfloor}}+{r_{ \lfloor\eta_{x, y}\psi_{x, y}^{(t)}\rfloor}})L_1  + f_{x}(y)-{f_{x^*}(y)}\right)^2}{2\left({\frac{(\sigma_{x^*}(y)+L_2c\lfloor\eta_{x^*, y}\psi_{x^*, y}^{(t)}\rfloor^{-\frac{\xi}{d}})^2}
			{\lfloor\eta_{x^* y}\psi_{x^*, y}^{(t)}\rfloor}
			+\frac{(\sigma_{x}(y)+L_2c\lfloor\eta_{x, y}\psi_{x, y}^{(t)}\rfloor^{-\frac{\xi}{d}})^2}
			{\lfloor\eta_{x, y}\psi_{x, y}^{(t)}\rfloor}}\right)})\}\Big]\\
	=&\liminf_{t\to\infty}\min_{x\neq x^*(y)}\min\Big[-\frac{1}{t^{\frac{1}{1+\xi}}}\ln [P\{ \Psi_{x^*, y}^{(t)}< \lfloor\eta_{x^*, y}\psi_{x^*, y} ^{(t)}\rfloor\}], -\frac{1}{t^{\frac{1}{1+\xi}}}\ln [P\{ \Psi_{x, y}^{(t)}<  \lfloor\eta_{x, y}\psi_{x, y} ^{(t)}\rfloor\}] , \\ 
	&\frac{1}{t^{\frac{1}{1+\xi}}}
	\frac{\left(-({r_{ \lfloor\eta_{x^*, y}\psi_{x^*, y}^{(t)}\rfloor}}+{r_{ \lfloor\eta_{x, y}\psi_{x, y}^{(t)}\rfloor}})L_1  + f_{x}(y)-{f_{x^*}(y)}\right)^2}{2\left({\frac{(\sigma_{x^*}(y)+L_2c\lfloor\eta_{x^*, y}\psi_{x^*, y}^{(t)}\rfloor^{-\frac{\xi}{d}})^2}
			{\lfloor\eta_{x^* y}\psi_{x^*, y}^{(t)}\rfloor}
			+\frac{(\sigma_{x}(y)+L_2c\lfloor\eta_{x, y}\psi_{x, y}^{(t)}\rfloor^{-\frac{\xi}{d}})^2}
			{\lfloor\eta_{x, y}\psi_{x, y}^{(t)}\rfloor}}\right)}\Big]\\
	=&C^\frac{1}{1+\xi}\min\{\min_{x\in\chi}\{\eta_{x,y}^{-\xi}(\eta_{x,y} ^{1+\xi}\ln\eta_{x,y} ^{1+\xi}-\eta_{x,y} ^{1+\xi}+ 1 )  [a(x,y)]^\frac{1}{1+\xi}\},\min_{x\neq x^*(y)}\{  \frac{(f_x(y)-f_{x^*}(y))^2}{2(\frac{\sigma_{x^*}^2(y)}{\eta_{x^* y}[a(x^*,y)]^\frac{1}{1+\xi}}+\frac{\sigma_{x}^2(y)}{\eta_{x, y}[a(x,y)]^\frac{1}{1+\xi}})}\}
\end{split}
\end{equation}
The first inequality holds due to (\ref{equation:CDLupper static})	. The first equality holds as $\max_{x\neq x^*(y)}d({f^\dagger}^{(t)}(y),{f^*(y)})$ is a globally bounded value and any bounded value will vanish under operations of $-\frac{1}{t^\frac{1}{1+\xi}}\ln[\cdot]$. The second equality holds by changing of the signature. The second inequality holds due to (\ref{inequality: insufficient samples}) and (\ref{inequation: mixed probability upper bound static}). The third equality holds as the constant $3$ will vanish under $-\frac{1}{t^\frac{1}{1+\xi}}\ln[\cdot]$. The last equality holds as each term in $\min$ converge to a constant and there are only finite number of terms. Notice that the bound holds for arbitrary $\eta_{x,y}\in(0,1)$, $\forall x\in\chi, y\in Y^\circ$. Therefore, to achieve the tightest bound:
\begin{equation}
\begin{split}
	&\liminf_{t\to\infty}-\frac{1}{t^{\frac{1}{1+\xi}}}\ln \upsilon^{(t)} (y)\\
	\geq&\sup_{\eta_{x,y}\in(0,1),\forall x\in\chi}\Big\{C^\frac{1}{1+\xi}\min\big\{\min_{x\in\chi}\{\eta_{x,y}^{-\xi}(\eta_{x,y} ^{1+\xi}\ln\eta_{x,y} ^{1+\xi}-\eta_{x,y} ^{1+\xi}+ 1 )  [a(x,y)]^\frac{1}{1+\xi}\},\\
	&\min_{x\neq x^*(y)}\{  \frac{(f_x(y)-f_{x^*}(y))^2}{2(\frac{\sigma_{x^*}^2(y)}{\eta_{x^* y}[a(x^*,y)]^\frac{1}{1+\xi}}+\frac{\sigma_{x}^2(y)}{\eta_{x, y}[a(x,y)]^\frac{1}{1+\xi}})}\}\big\}\Big\}\\
	\geq&\sup_{\eta_{x,y}\in(0,1),\forall x\in\chi}\Big\{C^\frac{1}{1+\xi}\min\big\{\min_{x\in\chi}\{\eta_{x,y}^{-\xi}(\eta_{x,y} ^{1+\xi}\ln\eta_{x,y} ^{1+\xi}-\eta_{x,y} ^{1+\xi}+ 1 )  [a(x,y)]^\frac{1}{1+\xi}\},\\
	&\min_{x\neq x^*(y)}\{  \frac{(f_x(y)-f_{x^*}(y))^2}{4\frac{\sigma_{x}^2(y)}{\eta_{x, y}[a(x,y)]^\frac{1}{1+\xi}}}, \frac{(f_x(y)-f_{x^*}(y))^2}{4\frac{\sigma_{x^*}^2(y)}{\eta_{x^*, y}[a(x^*,y)]^\frac{1}{1+\xi}}} \}  \big\}\Big\}\\			
	=&\sup_{\eta_{x,y}\in(0,1),\forall x\in\chi}\Big\{C^\frac{1}{1+\xi}\min\big\{\min_{x\in\chi}\{\eta_{x,y}^{-\xi}(\eta_{x,y} ^{1+\xi}\ln\eta_{x,y} ^{1+\xi}-\eta_{x,y} ^{1+\xi}+ 1 )  [a(x,y)]^\frac{1}{1+\xi}\},\\
	&\min_{x\neq x^*(y)}\{  \zeta_{x}(y)\eta_{x, y}[a(x,y)]^\frac{1}{1+\xi}\}, \zeta_{x^*,y}{\eta_{x^*, y}[a(x^*,y)]^\frac{1}{1+\xi}}\}  \big\}\Big\}\\			
	=&C^\frac{1}{1+\xi}\min_{x\in\chi}  \left\{ \zeta_{x}(y) \left(-W_{-1}\left(-\exp[{-    \zeta_{x}(y)   -1}]\right)\right)^{-\frac{1}{1+\xi}} [a(x,y)]^\frac{1}{1+\xi}\right\}    \\
	=:&H'_y(\{a(x,y)\}_{x\in\chi}),\\
\end{split}\label{definition: {H'_y} static}
\end{equation}
where
\begin{equation}
\zeta_{x}(y):=\begin{cases}
	\frac{(f_x(y)-f_{x^*}(y))^2}{4 \sigma_{x}^2(y)}& x\neq x^*(y)\\
	\frac{(   \min_{x\neq x^*(y)}  f_x(y)-f_{x^*}(y))^2}{4\sigma_{x^*}^2(y)} & x=x^*(y)
\end{cases}.
\end{equation}
The second equality holds as to achieve supremum, we must have $\eta_{x, y}$ satisfying $\eta_{x,y}^{-\xi}(\eta_{x,y} ^{1+\xi}\ln\eta_{x,y} ^{1+\xi}-\eta_{x,y} ^{1+\xi}+ 1 )  =\zeta_{x,y} \eta_{x, y}$ for each $x\in\chi$.
\endproof{$\blacksquare$}

\subsection{Proposition \ref{proposition: subset rate static}}
\proof{Proof}
The proof is similar to the proof of Proposition \ref{proposition: subset rate}. By substituting ${H_y}$ with $H'_y$, the statement is proved.
\endproof{$\blacksquare$}

\subsection{Theorem \ref{Theorem: optimal budget allocation}}
\proof{Proof}
The goal is to prove the optimality and feasibility of the following solution
\begin{equation}
a^*(x,y)=\frac{\lambda \left(-W_{-1}\left(-\exp[{-    \zeta_{x}(y)   -1}]\right)\right)}{C[\zeta_{x}(y)]^{1+\xi} }\quad\forall x\in\chi, y\in Y_{sub},
\end{equation}
where
\begin{equation}
\begin{split}
	\zeta_{x}(y):=&\begin{cases}
		\frac{(f_x(y)-f_{x^*}(y))^2}{4 \sigma_{x}^2(y)}& x\neq x^*(y)\\
		\frac{(   \min_{x\neq x^*(y)}  f_x(y)-f_{x^*}(y))^2}{4\sigma_{x^*}^2(y)} & x=x^*(y)
	\end{cases},\\
	\lambda=&\frac{1}{\sum_{x\in\chi}\int_{Y_{sub}}	\frac{\left(-W_{-1}\left(-\exp[{-    \zeta_{x}(y)   -1}]\right)\right)}{C [\zeta_{x}(y)]^{1+\xi} } dy}.
\end{split}
\end{equation}

For feasibility, $\forall x\in\chi$, as long as $\zeta_{x}(y)>0$, it is guaranteed that $a^*(x,y)>0$, $\forall y\in Y_{sub}$. We also have:
\begin{equation}
\sum_{x\in\chi}\int_{Y_{sub}} a^*(x,y)	dy=1
\end{equation}

For optimality, notice that when $a^*(x,y)$ is applied, for any $y$, $H'_y(\{a^*(x,y)\}_{x\in\chi})$ is the same:
\begin{equation}
\begin{split}
	&H'_y(\{a^*(x,y)\}_{x\in\chi})\\
	=&C^\frac{1}{1+\xi}\min_{x\in\chi}  \left\{ \zeta_{x}(y) \left(-W_{-1}\left(-\exp[{-    \zeta_{x}(y)   -1}]\right)\right)^{-\frac{1}{1+\xi}} [a^*(x,y)]^\frac{1}{1+\xi}\right\}\\
	=&\lambda^{1+\xi}.
\end{split}
\end{equation}
Thus, $H'_{Y_{sub}}(\{a^*(x,y)\}_{x\in\chi, y\in Y_{sub}})=\lambda^{1+\xi}$.

If there exists a better solution $a'(x,y)$ such that $H'_{Y_{sub}}(\{a'(x,y)\}_{x\in\chi, y\in Y_{sub}})>H'_{Y_{sub}}(\{a^*(x,y)\}_{x\in\chi, y\in Y_{sub}})$, we must have $a'(x,y)\geq a^*(x,y)$ at any $(x,y)$. Otherwise, if for some $x\in \chi, y\in Y_{sub}$, $a'(x,y)< a^*(x,y)$, then for this $y$:
\begin{equation}
\begin{split}
	&H'_y(\{a'(x,y)\}_{x\in\chi})\\
	=&C^\frac{1}{1+\xi}\min_{x\in\chi}  \left\{ \zeta_{x}(y) \left(-W_{-1}\left(-\exp[{-    \zeta_{x}(y)   -1}]\right)\right)^{-\frac{1}{1+\xi}} [a'(x,y)]^\frac{1}{1+\xi}\right\}\\
	<&C^\frac{1}{1+\xi}\min_{x\in\chi}  \left\{ \zeta_{x}(y) \left(-W_{-1}\left(-\exp[{-    \zeta_{x}(y)   -1}]\right)\right)^{-\frac{1}{1+\xi}} [a^*(x,y)]^\frac{1}{1+\xi}\right\}\\
	=&\lambda^{1+\xi},
\end{split}
\end{equation}
which suggests that 
\begin{equation}
H'_{Y_{sub}}(\{a'(x,y)\}_{x\in\chi, y\in Y_{sub}})\leq H'_y(\{a'(x,y)\}_{x\in\chi})<\lambda^{1+\xi}=H'_{Y_{sub}}(\{a^*(x,y)\}_{x\in\chi, y\in Y_{sub}}),
\end{equation}which is contradiction. Therefore, we must have  $a'(x,y)\geq a^*(x,y)$ for all $x\in \chi, y\in Y_{sub}$. However, to let $H'_{Y_{sub}}(\{a'(x,y)\}_{x\in\chi, y\in Y_{sub}})>H'_{Y_{sub}}(\{a^*(x,y)\}_{x\in\chi, y\in Y_{sub}})$, we must have $a'(x,y)> a^*(x,y)$ for all $x\in \chi, y\in Y_{sub}$. Otherwise, if for some $x\in \chi, y\in Y_{sub}$, $a'(x,y)= a^*(x,y)$, for this $x,y$:
\begin{equation}
\begin{split}
	&H'_{Y_{sub}}(\{a'(x,y)\}_{x\in\chi, y\in Y_{sub}})\\
	\leq&  H'_y(\{a'(x,y)\}_{x\in\chi})\\
	\leq& C^\frac{1}{1+\xi}      \zeta_{x}(y) \left(-W_{-1}\left(-\exp[{-    \zeta_{x}(y)   -1}]\right)\right)^{-\frac{1}{1+\xi}} [a'(x,y)]^\frac{1}{1+\xi}\\
	=&C^\frac{1}{1+\xi}      \zeta_{x}(y) \left(-W_{-1}\left(-\exp[{-    \zeta_{x}(y)   -1}]\right)\right)^{-\frac{1}{1+\xi}} [a^*(x,y)]^\frac{1}{1+\xi}\\
	=&\lambda^{1+\xi}\\
	=&H'_{Y_{sub}}(\{a^*(x,y)\}_{x\in\chi, y\in Y_{sub}}).
\end{split}
\end{equation}
Therefore, if such $a'(x,y)$ exists, we must have $a'(x,y)> a^*(x,y)$ for all $x\in \chi, y\in Y_{sub}$. However, all such allocations will violate the feasibility constraint:
\begin{equation}
\sum_{x\in\chi}\int_{Y_{sub}} a'(x,y)-a^*(x,y)dy>0. 
\end{equation}
The inequality holds as the integral of a strictly positive function on a domain with strictly positive measure is strictly positive. This leads to the fact that a better feasible $a'(x,y)$ does not exist. This suggests that $a^*(x,y)$ is the optimal.

\section{Algorithms}
\begin{algorithm} \label{alg: CONE}
\caption{CONE: Contextual Optimizer through Neighborhood Estimation}
\begin{algorithmic}[1]
	\State Set $\xi\in(0,1)$, and $c=\frac{\sup_{y_1,y_2\in Y^\circ}\{|y_1-y_2|\}}{2^{-\frac{\xi}{d}}}$ so that at any $x,y$, $B_2(x,y)\supset Y^\circ$ 
	\State Set $t=2|\chi|$, for each $x$, run simulation on two $y$ and get $\mathcal{D}^{(t)}$. Set a small constant $p_0\in(0,1)$
	\While{{\it Stopping criteria} isn't met:} 
	\While{ Not getting an {\it accepted} point} 
	\State Uniformly select one $(x,y)\in\chi\times Y^\circ$
	\State  At $(x,y)$, compute $\beta_{x, y}$ based on $\mathcal{D}^{(t)}$ \begin{align*}
		\beta_{x, y}=\begin{cases}
			\lambda_U&\frac{ -W_{-1}\left(-\exp[{-    \hat{\zeta}_{x}^{(t)}(y)   -1}]\right)}{[\hat{\zeta}_{x}^{(t)}(y)]^{1+\xi} }\geq 	\lambda_U\\
			\lambda_L&\frac{ -W_{-1}\left(-\exp[{-    \hat{\zeta}_{x}^{(t)}(y)   -1}]\right)}{[\hat{\zeta}_{x}^{(t)}(y)]^{1+\xi} }\leq 	\lambda_L\\
			\frac{ -W_{-1}\left(-\exp[{-    \hat{\zeta}_{x}^{(t)}(y)   -1}]\right)}{[\hat{\zeta}_{x}^{(t)}(y)]^{1+\xi} }&otherwise
		\end{cases}
	\end{align*}
	\State Rejection: uniformly generate $\tilde{\beta}\in[0, \lambda_U]$
	\If{$\beta_{x, y} < \tilde{\beta}$} 
	\State {\it Accept} $(x, y)$.
	\EndIf
	\EndWhile
	\State Run simulation on accepted $(x, y)$ and get $\tilde{G}_{[t+1]}$
	\State Update: $\mathcal{D}^{(t+1)}$, 
	\EndWhile
	\State Output ${\hat{x}^{(T)}(y)}$ at stopping stage $T$ based on $\mathcal{D}^{(T)}$ using SNE
\end{algorithmic}
\end{algorithm}

\begin{algorithm} \label{alg: USSNE}
\caption{US-SNE: Fully sequential Uniform Sampling for Shrinking Neighborhood Estimation}
\begin{algorithmic}[1]
	\State Set $\xi\in(0,1)$, and $c=\frac{\sup_{y_1,y_2\in Y^\circ}\{|y_1-y_2|\}}{2^{-\frac{\xi}{d}}}$ so that at any $x,y$, $B_2(x,y)\supset Y^\circ$ 
	\State Set $t=2|\chi|$, for each $x$, run simulation on two $y$ and get $\mathcal{D}^{(t)}$.
	\While{{\it Stopping criteria} isn't met:} 
	\State Uniformly sample $(x,y)\in\chi\times Y^\circ$
	\State Run simulation on $(x, y)$ and get $\tilde{G}_{[t+1]}$
	\State Update $\mathcal{D}^{(t+1)}$
	\EndWhile
	\State Output ${\hat{x}^{(T)}(y)}$ at stopping stage $T$ based on $\mathcal{D}^{(T)}$ using SNE
\end{algorithmic}
\end{algorithm}

\begin{algorithm} \label{alg: USKrig}
\caption{US-Krig: Uniform Sampling with batch replications for Kriging Estimation}
\begin{algorithmic}[1]
	\State Set batch size $l$
	\State Set $t=2|\chi|l$, for each $x$, run simulation replications on two $y$. Each $(x,y)$ will be simulated for $l$ replications. Get $\mathcal{D}^{(t)}$.
	\While{{\it Stopping criteria} isn't met:} 
	\State Uniformly sample one $(x,y)\in\chi\times Y^\circ$
	\State Run simulation on $(x, y)$ with $l$ replications and get $\tilde{G}_{[t+1]}$
	\State Update $\mathcal{D}^{(t+1)}$
	\EndWhile
	\State Output ${\hat{x}^{(T)}(y)}$ at stopping stage $T$ based on $\mathcal{D}^{(T)}$ using Kriging. The heteroscedastic noise at observed $(x,y)$ is estimated by the sample variance of replications at $(x,y)$. 
\end{algorithmic}
\end{algorithm}

\end{document}